%Authors: P.G. Casazza and N.J. Kalton

%Title: Uniqueness of unconditional bases in Banach spaces

%Filename: casazzakaltonunqunc.tex
%TeX: AMSTeX
%Length: 82217 bytes
%Received Date: 10/7/96
%SubjectClass: 46B15, 46B07
%Abstract: We prove a general result on complemented unconditional
%basic sequences
%in Banach lattices and apply it to give some new examples of
%spaces
%with unique unconditional basis.  We show that Tsirelson space
%and
%certain Nakano spaces have the unique unconditional bases.  We
%also construct an example of a space with a unique unconditional
%basis
%with a complemented subspace failing to have a unique
%unconditional
%basis.

%Citation: Preprint

%Special character check block
%32   space        33 ! exclam. pt.   34 " double quote  35 # sharp
%36 $ dollar       37 % percent       38 & ampersand     39 ' prime
%40 ( left paren.  41 ) rt. paren.    42 * asterisk      43 + plus
%44 , comma        45 - minus         46 . period        47 / divide
%58 : colon        59 ; semi-colon    60 < less than     61 = equal
%62 > greater than 63 ? question mark 64 @ at
%91 [ left bracket 92 \ backslash     93 ] right bracket 94 ^ caret
%95 _ underline    96 ` left single quote
%123 { left brace  124 | vertical bar 125 } right brace  126 ~ tilda

%Insert your TeX file starting here.

%%%%%%%%%%%%%%%%%%%%%%%%%%%%%%%%%%%%%%%%%%%%%%%%%%%%%%%%%%%%%%%%%%%%%

%%%
%%%           AMS-TeX 2.1
%%%
%%%
%%%%%%%%%%%%%%%%%%%%%%%%%%%%%%%%%%%%%%%%%%%%%%%%%%%%%%%%%%%%%%%%%%

%%%%%%%%%%%%%%%%%%%%%%%%%%%%%%%%%%%%%%%%%%%%%%%%%%%%%%%%%%%%%%%%%%%%%
%continuation of casazza.tex
\documentstyle{amsppt}
\magnification=\magstep1
\parindent=1em
\baselineskip 15pt
\hsize=12.3 cm
\vsize=18.5 cm
\NoRunningHeads
\pageno=1

\vsize=7.4in

%%%%%%%%%%%%%%%%%%%%%%%%%%%%%%%%%%%%%%%%%%%%%%%%%%%%%%%%%%%%%%%%%%%%%
%macros
\def\supp{\text{supp }}

\def\dim{\text{dim }}

\def\Ave{\mathop{\text {Ave}}}

%\def\Ave{\mathop{\bold E}}

 % triplenorm of #1
\def\colon{{:}\;}
%%%%%%%%%%%%%%%%%%%%%%%%%%%%%%%%%%%%%%%%%%%%%%%%%%%%%%%%%%%%%%%%%%%%%
%%%%%%%%%%%%%%%%%%%%%%%%%%%%%%%%%%%%%%%%%%%%%%%%%%%%%%%%%%%%%%%%%%%%%
\topmatter

\title Uniqueness of unconditional bases in Banach spaces
\endtitle
\author
P.G. Casazza
and N.J. Kalton
\endauthor
\address
Department of Mathematics,
University of Missouri,
Columbia, Mo.  65211, U.S.A.
\endaddress
\email pete\@math.missouri.edu, mathnjk\@mizzou1.missouri.edu
\endemail
\thanks Both authors were supported  by  NSF
Grant DMS-9201357 \endthanks
\subjclass
46B15, 46B07
\endsubjclass
\abstract

We  prove  a  general  result  on complemented unconditional
basic sequences in Banach lattices and apply it to give some
new examples of spaces with unique unconditional basis.
We show that Tsirelson space and certain Nakano spaces  have
unique unconditional bases.
We  also  construct  an  example  of  a  space with a unique
unconditional basis with a complemented subspace failing  to
have a unique unconditional basis.

\endabstract

\endtopmatter

\vskip10pt
\heading
1. Introduction
\endheading
\vskip10pt

A Banach space with an unconditional basis is said to have a
unique   unconditional   basis   if   any   two   normalized
unconditional bases are equivalent after a permutation.
It is  well-known that  $\ell_2$ has  a unique unconditional
basis (cf.  \cite{17}) and a classic result of Lindenstrauss
and Pe\l czynski \cite{18} asserts that the spaces  $\ell_1$
and  $c_0$  also  have  unique  unconditional  bases;  later
Lindenstrauss and Zippin \cite{21}  showed that this is  the
complete list of spaces  with symmetric bases for  which the
unconditional basis is unique.
Subsequently Edelstein and Wojtaszczyk \cite{10} showed that
direct sums  of $\ell_1,\ell_2$  and $c_0$  also have unique
unconditional bases.
In  1985,  Bourgain,  Casazza,  Lindenstrauss  and  Tzafriri
\cite{3} studied the classification problem for such spaces.
Their main results showed that  $\ell_1(\ell_2),c_0(\ell_1),
\ell_1(c_0),c_0(\ell_2)$   and   2-convexified    Tsirelson
$T^{(2)}$   have   unique   unconditional   bases  but  that
$\ell_2(\ell_1)$ and $\ell_2(c_0)$ do not.
Based  on  their  results  a  complete  classification looks
hopeless.
We also remark that a recent example of Gowers \cite{12} may
be easily shown to have unique unconditional basis.
Thus  there  are  many  ``pathological''  spaces with unique
unconditional basis.

In this paper  we will give  (Theorem 3.5) a  simple and, we
feel, useful characterization of complemented  unconditional
basic  sequences  in  Banach  sequence  spaces which are not
sufficiently   Euclidean   (i.e.   do   not  have  uniformly
complemented $\ell_2^n$'s).
This  theorem  is  the  discrete  analogue of Theorem 8.1 of
\cite{14}; in fact the  basic arguments are very  similar to
those given in \cite{16} and \cite{14}, but we have opted to
present a self-contained proof here.
We then use this result  and the recent work of  Wojtaszczyk
\cite{26}  to  give  some  more  examples  of fairly natural
spaces with unique unconditional basis.
In  Section  5,  we   introduce  the  class  of   left-  and
right-dominant bases and  use this notion  to show that  the
Nakano space $\ell(p_n)$ has a unique unconditional basis if
$p_n\downarrow  1$  and  $(p_{n}-p_{2n})\log  n$  is bounded
(there is a dual result if $p_n\uparrow \infty$).
We  also  show  that  Tsirelson  space  $T$  has  a   unique
unconditional basis (a question raised in \cite{3} p. 62).
In Section 6, we use similar techniques to show that certain
complemented subspaces of Orlicz sequence spaces have unique
unconditional bases.
Based on these examples we are able to resolve Problem  11.2
(p.104) of \cite{3}  by showing that  there is a  space with
unique  unconditional  basis  with  a  complemented subspace
(spanned  by  a  subsequence  of  the basis) failing to have
unique unconditional basis.

Also in Section 4, we use Theorem 3.5 to give a contribution
to  the  problem  of  uniqueness  of  unconditional bases in
finite-dimensional spaces.
Specifically,   we    prove   that    in   any    class   of
finite-dimensional  lattices  so  that  $\ell_2^n$  is   not
complementably    and    disjointly    representable,    the
unconditional basis  is almost  unique; for  a more  precise
statement see Theorem 4.1.

We remark that the  techniques developed here using  Theorem
3.5  can  be  used  successfully  to obtain other results on
uniqueness.
In  particular  we  plan  to  study  unconditional  bases in
$c_0-$products in a later publication.
Since the  arguments in  such spaces  are considerably  more
complicated,  it  seemed,  however,  appropriate to restrict
attention here to some simple applications.

\vskip10pt \heading{2.  Definitions and notation}\endheading
\vskip10pt

We will take the viewpoint that an unconditional basis in  a
Banach space $X$ confers  the structure of an  atomic Banach
lattice on $X$.
We  will  thus  adopt  the  language and structure of Banach
lattices.
It  is  well-known  that  a  separable Banach lattice can be
regarded as a K\"othe function space.

We will in general use the same notation as in \cite{16}.
Let $\Omega$ be  a Polish space  (i.e. a separable  complete
metric  space)  and  let  $\mu$  be  a $\sigma-$finite Borel
measure on $\Omega.$
We denote by  $L_0(\mu)$ the space  of all Borel  measurable
functions on $\Omega$, where we identify functions differing
only on a set of measure zero; the natural topology of $L_0$
is convergence in measure on sets of finite measure.
An admissible
norm is then a lower-semi-continuous map $f\to \|f\|$ from
$L_0(\mu)$ to $[0,\infty]$ such that:\newline \vskip.2truecm
\flushpar (a) $\|\alpha f\| =|\alpha|\|f\|$ whenever
$\alpha\in\bold R,\ f\in L_0.$ \newline (b) $\|f+g\|\le
\|f\|+\|g\|,$ for $f,g\in L_0.$ \newline (c) $\|f\| \le
\|g\|,$ whenever $|f| \le |g|$ a.e.  (almost
everywhere).\newline (d) $\|f\|<\infty$ for a dense set of
$f\in L_0,$ \newline (e) $\|f\|=0$ if and only if $f=0$ a.e.
\newline \vskip.1truecm

A K\"othe function space on $(\Omega,\mu)$ is defined
to be a dense order-ideal $X$ in $L_0(\mu)$ with an
associated admissible norm $\|\,\|_X$ such that if
$X_{max}=\{f:\|f\|_X<\infty\}$ then either:  \newline (1)
$X=X_{max}$ ($X$ is {\bf  maximal}) or:\newline (2) $X$ is
the closure of the simple functions in $X_{max}$ ($X$ is
{\bf  minimal}).\newline
Any order-continuous K\"othe function space is
minimal.
Also any K\"othe function space which does
not contain a copy of $c_0$ is both maximal and minimal.

If $X$ is an order-continuous K\"othe function space then
$X^*$ can be identified with the K\"othe function space of
all $f$ such that:  $$ \|f\|_{X^*}=\sup_{\|g\|_X\le 1} \int
|fg|\,d\mu <\infty.$$ $X^*$ is always maximal.

A K\"othe function space $X$ is said to be $p-$convex (where
$1<p<\infty)$ if there is a  constant $C$ such that for  any
$f_1,\ldots,f_n\in X$ we have
$$
\|(\sum_{i=1}^n|f_i|^p)^{1/p}\|_X                        \le
C(\sum_{i=1}^n\|f_i\|_X^p)^{1/p}.
$$
$X$ is said  to have an  upper $p$-estimate if  for some $C$
and any disjoint $f_1,\ldots,f_n\in X,$
$$
\|\sum_{i=1}^n f_i\|_X \le C(\sum_{i=1}^n\|f_i\|_X^p)^{1/p}.
$$
$X$ is  said to  be $q-$concave  ($0<q<\infty$) if  for some
$c>0$ and any $f_1,\ldots,f_n\in X$ we have
$$
\|(\sum_{i=1}^n|f_i|^q)^{1/q}\|_X                        \ge
c(\sum_{i=1}^n\|f_i\|_X^q)^{1/q}.
$$
$X$ is said to have  a lower $q$-estimate if for  some $c>0$
and any disjoint $f_1,\ldots,f_n\in X,$
$$
\|\sum_{i=1}^n       f_i\|_X       \ge        c(\sum_{i=1}^n
\|f_i\|_X^q)^{1/q}.
$$

A Banach space  $X$ is said  to be of  (Rademacher) type $p$
($1\le p\le 2$) if there is  a constant $C$ so that for  any
$x_1,\ldots,x_n \in X,$
$$
\mathop{\text{Ave    }}_{\epsilon_i=\pm     1}\|\sum_{i=1}^n
\epsilon_i x_i\| \le C(\sum_{i=1}^n\|x_i\|^p)^{1/p}
$$
and $X$ is of cotype $q$ ($2\le q<\infty$) if for some $c>0$
and any $x_1,\ldots,x_n\in X$ we have
$$
\mathop{\text{Ave
}}_{\epsilon_i=\pm1}\|\sum_{i=1}^n\epsilon_ix_i\|        \ge
c(\sum_{i=1}^n\|x_i\|^q)^{1/q}.
$$
We recall that a Banach lattice has nontrivial cotype  (i.e.
has cotype $q<\infty$  for some $q$)  if and only  if it has
nontrivial   concavity   (i.e.   is   $q-$concave  for  some
$q<\infty$).
If $X$ is  a Banach lattice  which has nontrivial  concavity
then  there  is  a   constant  $C=C(X)$  so  that   for  any
$x_1,\ldots,x_n\in X$ we have
$$
\frac1C
(\Ave_{\epsilon_k=\pm1}\|\sum_{k=1}^n\epsilon_kx_k\|^2)^{1/2}\le
\|(\sum_{k=1}^n|x_k|^2)^{1/2}\|_X           \le            C
(\Ave_{\epsilon_k=\pm1}\|\sum_{k=1}^n\epsilon_kx_k\|^2)^{1/2}.
$$

We  will  use  the  term  {\bf   sequence  space} to denote a
K\"othe  function  space  $X$  on  $\bold  N$  equipped with
counting   measure,   and   subject   to  the  normalization
constraint that  $\|e_j\|_X=1$ for  all $j\in\bold  N$ where
$e_j=\chi_{\{j\}}.$
It is clear that $(e_n)$ forms an unconditional basis for  a
sequence  space  $X$  if  and  only  if  $X$  is minimal (or
separable).
We will consider finite-dimensional sequence spaces modelled
on finite sets $[N]=\{1,2,\ldots,N\}$ with counting measure.

In keeping with current usage we will write $c_{00}$ for the
space of finitely nonzero sequences.
If $A$ is a subset of  $\bold N$ we write $e_A$ in  place of
$\chi_A$ and if $x$ is  any sequence we write $Ax=e_Ax.$  If
$A,B$  are  subsets  of  $\bold  N$  we write $A<B$ if $a<b$
whenever $a\in A$ and $b\in B.$
If $x$ is a sequence then $\supp x=\{i\colon x(i)\neq 0\}.$

Many of our examples will be Orlicz sequence spaces or  more
general Orlicz-Musielak or modular sequence spaces.
If  $(F_n)$  is  a  sequence  of  Orlicz  functions then the
modular  sequence  space  $\ell_{(F_n)}$  is  the  space  of
sequences       $(x(n))_{n=1}^{\infty}$       such      that
$\sum_{n=1}^{\infty}F_n(|x(n)|)<\infty,$ with the norm
$$
\|x\|_{\ell(F_n)}=\inf\{\lambda>0\colon
\sum_{n=1}^{\infty}F_n(\lambda^{-1}|x(n)|)\le 1\}.
$$
In the  case $F_n=F$  for all  $n$ we  have the Orlicz space
$\ell_F.$
If $\ell_{(F_n)}$ is separable or has finite cotype then the
canonical  basis  vectors  form  an  unconditional  basis of
$\ell_{(F_n)}$; otherwise they  form an unconditional  basis
of their closed linear span $h_{(F_n)}.$
We refer to \cite {19}  for the basic properties of  modular
sequence spaces.

One special  case is  to take  $F_n(t)=t^{p_n}$ where  $1\le
p_n<\infty.$
This  is  often  called  a  Nakano  space  and  we denote it
$\ell(p_n).$
$\ell(p_n)$ is separable if and only if $\sup p_n<\infty.$
It may also be shown that  if $p_n>1$ for all $n$ and  $\sup
p_n<\infty$     then     $\ell{(p_n)}^*=\ell(q_n)$     where
$p_n^{-1}+q_n^{-1}=1.$
If $\sup p_n=\infty$ then  we write $h(p_n)$ for  the closed
linear   span   of   the   basis   vectors,   and   we  have
$h(p_n)^*=\ell(q_n).$

Let  $(u_n)$   and  $(v_n)$   be  two   unconditional  basic
sequences.
We  say  that  $(u_n)$  and  $(v_n)$  are {\bf  permutatively
equivalent} if there is a permutation $\pi$ of $\bold N$  so
that $(u_n)$ and $(v_{\pi(n)})$ are equivalent.
We say  that $(u_n)$  is {\bf   equivalent to  its square} if
$(u_n)$   is   permutatively   equivalent   to   the   basis
$\{(u_1,0),(0,u_1),(u_2,0),\ldots\}$ of $[u_n]\oplus[u_n].$
A Banach space  $X$ with an  unconditional basis has  {\bf  a
unique   unconditional   basis}   if   any   two  normalized
unconditional bases are permutatively equivalent.
We remark that there  is an important Cantor-Bernstein  type
principle which  helps determine  whether two  unconditional
bases  are   permutatively  equivalent:     if   $(u_n)$  is
permutatively equivalent  to some  subset of  $(v_n)$ and if
$(v_n)$  is  permutatively  equivalent  to  some  subset  of
$(u_n)$   then   $(u_n)$   and   $(v_n)$  are  permutatively
equivalent.
We are grateful to P. Wojtaszczyk for drawing our  attention
to this principle, which appears explicitly in \cite{27} and
is used in \cite{26}.
We are indebted to C.  Bessaga for the information that  the
Cantor-Bernstein principle  was used  implicitly earlier  by
Mityagin in \cite{22}.

A Banach space $X$ is called {\bf  sufficiently Euclidean} if
there  is  a  constant  $M$  so  that  for any $n$ there are
operators $S\colon X\to\ell_2^n$ and $T\colon \ell_2^n\to X$
so that $ST=I_{\ell_2^n}$ and $\|S\|\|T\|\le M.$
We will say  that $X$ is  {\bf  anti-Euclidean} if  it is not
sufficiently Euclidean.

A Banach  lattice $X$  is called  {\bf sufficiently  lattice
Euclidean} if there  is a constant  $M$ so that  for any $n$
there    are    operators    $S\colon    X\to\ell_2^n$   and
$T\colon\ell_2^n\to   X$   so   that  $ST=I_{\ell_2^n}$  and
$\|S\|\|T\|\le  M,$   and  such   that  $S$   is  a  lattice
homomorphism.
This  is  equivalent  to  asking  that  $\ell_2$ is finitely
representable as a complemented sublattice of $X.$
We will say that $X$  is {\bf lattice anti-Euclidean} if  it
is not sufficiently lattice Euclidean.
We  use  the  same  terminology  for  an unconditional basic
sequence, which we regard as inducing a lattice structure on
its closed linear span.

Finally if $\Cal  X$ is a  family of Banach  lattices we say
that $\Cal  X$ is  {\bf sufficiently  lattice Euclidean}  if
there is  a constant  $M$ so  that for  any $n$ there exists
$X\in\Cal  X$  and  operators  $S\colon  X\to\ell_2^n$   and
$T\colon  \ell_2^n\to  X$  so  that  $ST=I_{\ell_2^n}$   and
$\|S\|\|T\|\le  M,$   and  such   that  $S$   is  a  lattice
homomorphism.
If $\Cal X$  is not sufficiently  lattice Euclidean we  will
say that it is {\bf lattice anti-Euclidean}.

\vskip1truein

\heading
3. Complemented unconditional basic sequences
\endheading

The main results of this  section are Theorems 3.4 and  3.5,
which   show   that   complemented   lattice  anti-Euclidean
unconditional basic sequences in an order-continuous  Banach
lattice or Banach sequence space take a particularly  simple
form.

\proclaim{Lemma 3.1} Let $X$ be a Banach sequence space  and
suppose $(u_1,\ldots,u_n)$ are  disjoint elements of  $X_+$,
and $(u_1^*,\ldots,u_n^*)$ are disjoint in $X^*_+$.
Suppose that $M\ge 1$ is a constant such that
$$
\|\sum_{j=1}^na_ju_j\|_X \le M(\sum_{j=1}^n|a_j|^2)^{1/2}
$$
and
$$
\|\sum_{j=1}^na_ju_j^*\|_{X^*}                           \le
M(\sum_{j=1}^n|a_j|^2)^{1/2}
$$
whenever $a_1,\ldots,a_n\in\bold R.$
Suppose further that
$$
\sum_{j=1}^n\langle u_j,u_j^*\rangle =\alpha n.
$$
Then  $\ell_2^m$  is  $2M^2\alpha^{-1}$-representable  as  a
$2M^2\alpha^{-1}$-complemented sublattice  of $X,$  for some
$m\ge \frac12\alpha M^{-2}n.$
\endproclaim

\demo{Proof}We can clearly suppose that $\supp u=\supp u^*.$
Note that $\langle u_i,u_i^*\rangle \le M^2.$
Let      $J=\{j\colon      \langle       u_j,u_j^*\rangle\ge
\frac12\alpha\}.$
Then $|J|\ge\frac12 M^{-2} \alpha n.$
Notice that for any $(a_j)_{j\in J}$ we have
$$
\|\sum_{j\in  J}a_ju_j\|_X\|  \sum_{j\in  J}a_ju_j^*\|_{X^*}
\ge \frac{\alpha}{2}(\sum_{j\in J}|a_j|^2).
$$
Thus
$$
\frac{\alpha}{2M}(\sum_{j\in      J}|a_j|^2)^{1/2}       \le
\|\sum_{j\in J} a_ju_j\|_X \le M(\sum_{j\in J}|a_j|^2)^{1/2}
$$
so   that   $[u_j]_{j\in   J}$   is  $(2M^2)/\alpha$-lattice
isomorphic to $\ell_2^{|J|}.$
If we let  $\gamma_j=\langle u_j,u^*_j\rangle^{-1}$ then  we
can define a projection $P$ onto $[u_j]_{j\in J}$ by
$$
Px=\sum_{j\in J}\gamma_j^{-1}\langle x,u_j^*\rangle u_j.
$$
Then
$$
\|Px\|_X    \le    \frac{2M}{\alpha}(\sum_{j\in   J}|\langle
x,u_j^*\rangle|^2)^{1/2}.
$$
However,
$$
\align  \sum_{j\in  J}|\langle  x,u_j^*\rangle|^2  &=\langle
x,\sum_{j\in J}\langle  x,u_j^*\rangle u_j^*\rangle  \\ &\le
M\|x\|_X   (\sum_{j\in   J}\langle   x,u_j^*\rangle^2)^{1/2}
\endalign
$$
whence we obtain
$$
\|Px\|_X \le \frac{2M^2}{\alpha}\|x\|_X.
$$
\qed\enddemo

\proclaim{Lemma 3.2}Let  $X$ be  a order-continuous  K\"othe
function space  on $(\Omega,\mu)$.   Suppose  $m\in\bold N,$
and $\phi\in L_1(\mu)$ with $\phi\ge 0.$
Suppose      $f_1,\ldots,f_n\in      X_+$      and       let
$F=(\sum_{j=1}^nf_j^2)^{1/2}$ and $F_{\infty}=\max_jf_j.$
Then we can  partition $[n]=\{1,2,\ldots,n\}$ into  $m$-sets
$J_1,\ldots,J_m$  and  find  a  set  $A\subset  \Omega$ with
$\int_A\phi\,d\mu\ge \frac34\int_{\Omega}\phi\,d\mu$ so that
whenever $a_1,\ldots,a_m\in\bold R$ we have
$$
\|(\sum_{k=1}^m a_k^2\sum_{j\in J_k}f_j^2)^{1/2}\chi_A  \|_X
\le            2(\|F\|_X+             2.5^m\|F_{\infty}\|_X)
m^{-1/2}(\sum_{k=1}^ma_k^2)^{1/2}.
$$
\endproclaim

\demo{Proof}
We    may    select    a    collection   of   $5^m$   points
$(b^r)_{r=1}^{5^m}$ in the unit sphere of $\ell_2^m$ to form
a $\frac12$-net.
Then  if  $T\colon\ell_2^m\to  Y$  is  any  operator we have
$\|T\|\le 2\sup_{r\le 5^m}\|Tb^r\|_Y.$

Next let $\Pi$ be the  set of all $m^n$ partitions  of $[n]$
into  $m$-sets  $\pi=(J_1,\ldots,J_m)$  with  a  probability
measure $P$ defined to be normalized counting measure.
For   $1\le   l\le   m$   and   $1\le   j\le  n$  we  define
$\xi_{lj}(\pi)=1$ if $j\in J_l$ and $0$ otherwise.
The random variables $\xi_{ki},\xi_{lj}$ are independent  if
$i\neq j$,  and each  have expectation  $1/m,$ and  variance
$(m-1)/m^2.$
We will use the fact  that the covariance of $\xi_{li}$  and
$\xi_{ki}$ is negative:  it can be computed as $-1/m^2.$

For   each   $\pi\in\Pi$   we   set    $x_l(\pi)=(\sum_{j\in
J_l}f^2_j)^{1/2}.$  Then  for  $1\le  r\le  5^m$  we  define
$A(\pi,r)$ to be  the set of  $\omega\in\Omega$ such that  $
(\sum_{l=1}^m|b^r_l|^2|x_l(\pi,\omega)|^2)^{1/2}         \le
m^{-1/2}(F(\omega)+ 2.5^mF_{\infty}(\omega)).$

Let us fix  $\omega$ and $1\le  r\le 5^m,$ and  consider the
random variable $\zeta$ on $\Pi$ defined by
$$
\zeta(\pi) = \sum_{l=1}^m |b^r_l|^2\sum_{j=1}^n\xi_{lj}(\pi)
|f_j(\omega)|^2= \sum_{l=1}^m|b^r_l|^2|x_l(\pi,\omega)|^2.
$$
Then
$$
\int \zeta \,dP =\frac1m F(\omega)^2.
$$
Next  we  estimate  the  variance  of  $\zeta$ recalling our
previous  observations   concerning  the   random  variables
$\xi_{ki}.$
$$
\int     (\zeta-\frac1mF(\omega)^2)^2\,dP     \le    \frac1m
\sum_{l=1}^m|b^r_l|^4\sum_{j=1}^n|f_j(\omega)|^4 \le \frac1m
F(\omega)^2F_{\infty}(\omega)^2.
$$

Now
$$
\int_{\omega\notin   A(\pi,r)}(\zeta-\frac1mF(\omega)^2)^2dP
\ge         \frac{4.5^m}{m}F(\omega)^2F_{\infty}(\omega)^2P(
\omega\notin A(\pi,r)).
$$
Thus
$$
P(\omega\notin A(\pi,r)) \le \frac{1}{4.5^m},
$$
Let $B(\pi,r)$ be the complement of $A(\pi,r).$
Then
$$
\frac1{m^n}\sum_{\pi\in\Pi}
\int_{\Omega}\phi\chi_{B(\pi,r)}\,d\mu\le
\frac1{4.5^m}\int_{\Omega}\phi\,d\mu.
$$

Summing over  $r$ we  further obtain  the existence  of some
$\pi$ so that
$$
\int_{\Omega}\phi\sum_{r=1}^{5^m}\chi_{B(\pi,r)}\,d\mu   \le
\frac{1}{4}\int_{\Omega}\phi\,d\mu.
$$

For this fixed $\pi$, let $B=\cup_{r=1}^{5^m}B(\pi,r).$
Then   $\int_B\phi\,d\mu   \le   \frac14\int_{\Omega}   \phi
\,d\mu.$
Let $A$ be the complement of $B$.
Then $\int_{A}\phi\,d\mu\ge \frac34\int_{\Omega}\phi\,d\mu.$

For $1\le r\le 5^m$ we then have that
$$
\chi_A(\sum_{l=1}^{5^m}|b^r_l|^2x_l(\pi)^2)^{1/2}\le
m^{-1/2}(F+2.5^mF_{\infty}).
$$

By  considering   the  map   $T\colon\ell_2^m\to  X(\ell_2)$
defined by $T(e_l)=(0,\ldots,0,x_l,0,\ldots)$ with $x_l$  in
the   $l$th.   position   it   follows   that   for    every
$a_1,\ldots,a_m$ we have
$$
\|(\sum_{l=1}^m|a_l|^2|x_l|^2)^{1/2}\chi_A\|_X           \le
2m^{-1/2}(\|F\|_X+2.5^m\|F_{\infty}\|_X)
(\sum_{l=1}^m|a_l|^2)^{1/2}.\qed
$$
\enddemo

\proclaim{Lemma 3.3}Let $\Cal X$ be a lattice anti-Euclidean
family of Banach sequence spaces.
Then, given $M$  there exists $\delta=\delta(M)>0$  with the
following property.
Suppose  that  $Y$  is  an order-continuous K\"othe function
space on $(\Omega,\mu)$ and $X\in\Cal X.$
Let $d=\dim X\le \infty.$
Suppose $S\colon X\to  Y$ and $T\colon  Y\to X$ are  bounded
operators with $\|S\|,\|T\|\le M.$
Suppose  $(h_n)_{n=1}^d\in  L_1(\mu)$  satisfy  $0\le h_n\le
|Se_n||T^*e_n|$ and $\int h_n\,d\mu\ge M^{-1}.$
Then, for any $N\le d,$ and any $\alpha_1,\ldots,\alpha_N,$
$$
\int_{\Omega}\max_{1\le   k\le   N}|\alpha_kh_k|\,d\mu   \ge
\delta\sum_{k=1}^N|\alpha_k|.
$$
\endproclaim

\demo{Proof}
Let us suppose that for some $M$ the conclusion of the lemma
is false.
Suppose $m\in\bold N$ is given.
We put $\epsilon=(2.5^m)^{-1}$.
Then we can find $X\in\Cal X$ and $S\colon X\to Y,\  T\colon
Y\to   X$   with   $\|S\|,\|T\|\le   M$   and  $0\le  h_n\le
|Se_n||T^*e_n|$  for  $1\le   n\le  d=\dim  X$   with  $\int
h_nd\mu\ge M^{-1}$ and such that for suitable  $0\le\zeta\in
c_{00}$  and   $N\in\bold  N,$   with  $\|\zeta\|_1=1$   and
$\zeta(i)=0$ for all $i> N,$ we have
$$
\int_{\Omega}          \max_{1\le          k\le           N}
\zeta(k)h_k\,d\mu<\frac{\epsilon^2}{16M^4}.
$$

Let $f_n=Se_n$ and $g_n=T^*e_n.$
Notice that this implies that
$$
\|h_n\|_1\le \|Se_n\|_Y\|T^*e_n\|_{Y^*}\le M^2.
$$

It follows from Krivine's theorem (\cite{20} Theorem  1.f.4,
p.93) that if $\alpha_1,\ldots,\alpha_n\in\bold R$ then
$$
\|(\sum_{k=1}^n\alpha_k^2f_k^2)^{1/2}\|_Y                \le
K_GM\|\sum_{k=1}^n\alpha_ke_k\|_X,
$$
$$
\|(\sum_{k=1}^n\alpha_k^2g_k^2)^{1/2}\|_{Y^*}            \le
K_GM\|\sum_{k=1}^n\alpha_ke_k\|_{X^*},
$$
where $K_G$ is as usual the Grothendieck constant.

By  a  well-known  theorem  of  Lozanovskii we can factorize
$\zeta=\xi\xi^*$ where $0\le  \xi,\xi^*\in c_{00}$ have  the
same      support      as      $\zeta$      and      satisfy
$\|\xi\|_X=\|\xi^*\|_{X^*}=1.$

Next let $F=(\sum_{k=1}^N\xi(k)^2f_k^2)^{1/2}\in Y.$
It follows from the remarks above that $\|F\|_Y\le K_GM.$
Similarly if $G=(\sum_{k=1}^N\xi^*(k)^2g_k^2)^{1/2}\in  Y^*$
then $\|G\|_{Y^*}\le K_GM.$
Finally let $H=\sum_{k=1}^n\zeta(k)h_k.$
Then $M^{-1}\le \|H\|_1=\int H\,d\mu\le M^2.$

For  each  $k$  let  $B_k=\{  \xi(k)f_k\ge \epsilon F\}$ and
$B_k^*=\{\xi^*(k)g_k\ge \epsilon G\}.$
We will let $A_k=\Omega\setminus (B_k\cup B_k^*).$

Now  if  $\omega\in  B_k$  we  have $F\xi^*(k)g_k(\omega)\le
\epsilon^{-1}\zeta(k)h_k(\omega)$ and so
$$
\align        \sum_{k=1}^N\zeta(k)h_k\chi_{B_k}         &\le
(\sum_{k=1}^N\xi(k)^2f_k^2\chi_{B_k})^{1/2}(\sum_{k=1}^N\xi^*
(k)^2g_k^2)^{1/2}                  \\                   &\le
F(\sum_{k=1}^N\xi^*(k)^2g_k^2)^{1/2}\\                  &\le
\epsilon^{-1}(\sum_{k=1}^N\zeta(k)^2h_k^2)^{1/2}\\      &\le
\epsilon^{-1}H^{1/2}(\max_{1\le  k\le  N}\zeta(k)h_k)^{1/2}.
\endalign
$$
From this we deduce that
$$
\|\sum_{k=1}^N\zeta(k)h_k\chi_{B_k}\|_1   \le   \frac1{4M^2}
\|H\|_1^{1/2}\le \frac1{4M}.
$$

Similarly
$$
\|\sum_{k=1}^N\zeta(k)h_k\chi_{B^*_k}\|_1 \le \frac1{4M}.
$$

Hence
$$
\int \sum_{k=1}^N\zeta(k)h_k\chi_{A_k}d\mu \ge \frac1{2M}.
$$

Now $\max_{1\le j\le N}\xi(j)f_j\chi_{A_j}\le \epsilon F$ while
$\max_{1\le
j\le N}\xi^*(j)g_j\chi_{A_j}\le \epsilon G.$
Consider $X\oplus X^*$ (with the maximum norm) as a  K\"othe
function space on two copies of $(\Omega,\mu)$ and  consider
the functions  $(\xi(j)f_j\chi_{A_j},\xi^*(j)g_j\chi_{A_j})$
in $X\oplus X^*$.
Using              Lemma              3.2               with
$\phi=(\sum_{k=1}^N\zeta(k)h_k\chi_{A_k},\sum_{k=1}^Nh_k\chi_{A_k})$
it is easy to deduce the existence of a Borel subset $D$  of
$\Omega$ with
$$
\int_D \sum_{k=1}^N \zeta(k)h_k\chi_{A_k}d\mu \ge \frac1{2M}
$$
and a partition  $J_1,\ldots,J_m$ of $[N]$  so that for  any
$a_1,\ldots,a_m$ we have
$$
\|(\sum_{k=1}^ma_k^2\sum_{j\in
J_k}\xi(j)^2f_j^2\chi_{A_j})^{1/2}\chi_D\|_X             \le
4K_GMm^{-1/2}(\sum_{k=1}^ma_k^2)^{1/2}
$$
and
$$
\|(\sum_{k=1}^ma_k^2\sum_{j\in
J_k}\xi^*(j)^2g_j^2\chi_{A_j})^{1/2}\chi_D\|_{X^*}       \le
4K_GMm^{-1/2}(\sum_{k=1}^ma_k^2)^{1/2}.
$$

Let   $L$   be   the   set   of   $1\le   j\le  N$  so  that
$\|h_j\chi_{A_j\cap D}\|_1\ge\frac1{4M}.$
Then
$$
\int_D    \sum_{j\notin    L}\zeta(j)h_j\chi_{A_j}d\mu   \le
\frac1{4M}
$$
so that
$$
\int_D\sum_{j\in L}\zeta(j)h_j\chi_{A_j}d\mu \ge \frac1{4M}.
$$

Now let
$$
u_k=m^{1/2}\sum_{j\in J_k\cap L}\xi(j)e_j\in X
$$
and
$$
u_k^*=m^{1/2}\sum_{j\in J_k\cap L}\xi^*(j)e_j\in X^*.
$$

Consider  an  element  $v=\sum_{k=1}^ma_ku_k\in  X$  and let
$v^*=\sum_{j\in   L}v^*(j)e_j\in   X^*$   norm   $v$,   i.e.
$\|v^*\|_{X^*}=1$ and $\langle v,v^*\rangle=\|v\|_X.$

Then
$$
\align    \|v\|_X    &\le    4M\sum_{j\in   L}v(j)v^*(j)\int
h_j\chi_{A_j\cap    D}d\mu    \\    &\le   4M\int(\sum_{j\in
L}v(j)^2f_j^2\chi_{A_j\cap              D})^{1/2}(\sum_{j\in
L}v^*(j)^2g_j^2)^{1/2}d\mu\\      &\le       4M\|(\sum_{j\in
L}v(j)^2f_j^2\chi_{A_j\cap    D})^{1/2}\|_X    \|(\sum_{j\in
L}v^*(j)^2g_j^2)^{1/2}\|_{X^*}.  \endalign
$$

Here    the    first    factor    can    be   estimated   by
$4K_GM(\sum_{k=1}^Na_k^2)^{1/2}$  and  the  second factor by
Krivine's theorem is majorized by $K_GM.$
Hence
$$
\|v\|_X \le 2^4K_G^2M^3(\sum_{k=1}^ma_k^2)^{1/2}.
$$
Thus we have the inequality ' $$
\|\sum_{k=1}^ma_ku_k\|_X                                 \le
2^4K_G^2M^3(\sum_{k=1}^ma_k^2)^{1/2}.
$$

Precisely    dual    arguments    will    yield    that   $$
\|\sum_{k=1}^ma_ku_k^*\|_{X^*}\le
2^4K_G^2M^2(\sum_{k=1}^ma_k^2)^{1/2}.  $$

Finally $\sum_{k=1}^m\langle u_k,u_k^* \rangle  =m\sum_{j\in
L}\zeta(j)\ge     \frac14m\int_D\sum_{j\in     L}\zeta(j)h_j
\chi_{A_j} \ge 2^{-4}M^{-1}m.$

Now    Lemma    3.1    yields    that    $\ell_2^{m_0}$   is
$2^9K_G^4M^7$-complementably $2^9K_G^4M^7$-lattice  finitely
representable      in      $X$      for     some     $m_0\ge
2^{-9}K_G^{-4}M^{-5}m.$

It  is  clear  that  this  impossible  for arbitrarily large
$m.$\qed\enddemo

\proclaim{Theorem    3.4}Let    $Y$    be    a     nonatomic
order-continuous Banach lattice and suppose that $(f_n)$  is
a complemented unconditional basic sequence in $Y.$
Suppose $(f_n)$ is lattice anti-Euclidean.
Then  $(f_n)$  is  equivalent  to  a  complemented  disjoint
sequence $(f'_n)$ in $Y$.
\endproclaim

\demo{Proof}We  suppose  that  $Y$  is  an  order-continuous
K\"othe  function  space  on  $(\Omega,\mu),$ where $\mu$ is
nonatomic.
Let $X$ be  the sequence space  induced by $(f_n),$  and let
$S\colon X\to Y$ be the bounded linear map with $Se_n=f_n.$
Then there  is also  a bounded  linear map  $T\colon Y\to X$
with $TS=I_X.$
As before let $g_n=T^*e_n$ and $h_n=|f_ng_n|.$
Then for suitable $\delta>0$ we have
$$
\int_{\Omega}\max_{1\le     k\le     N}|\alpha_kh_k|     \ge
\delta\sum_{k=1}^N|\alpha_k|
$$
for every $N, \alpha_1,\ldots,\alpha_N.$
By a result of Dor \cite{9} there exist disjoint Borel  sets
$(E_n)_{n=1}^{\infty}$ so that $\int_{E_n}h_nd\mu= \delta.$
It      is      then      easy      to      verify      that
$(|f_n|\chi_{E_n})_{n=1}^{\infty}$    is    a   complemented
disjoint sequence equivalent to $(f_n)$.
Indeed define $U\colon  X\to Y$ by  $Ue_n=f_n\chi_{E_n}$ and
$V\colon  Y   \to  X$   by  $V(y)(j)   =  \delta^{-1}\langle
y,|g_j|\chi_{E_j}\rangle.$
Then $VU=I_X$ and for any $\xi\in c_{00}$ we have
$$
\|U\xi\|_Y                                               \le
\|(\sum_{j=1}^{\infty}\xi(j)^2f_j^2)^{1/2}\|_Y           \le
K_GM\|\xi\|_X.
$$
Also if $\xi^* \in c_{00}$ then
$$
|\langle       Vy,\xi^*\rangle|\le        \delta^{-1}\langle
|y|,(\sum_{j=1}^{\infty}\xi^*(j)^2g_j^2)^{1/2}\rangle    \le
\delta^{-1}K_GM\|\xi^*\|_{X^*}.
$$
Thus $U,V$  are both  bounded operators  and the  theorem is
proved.\qed\enddemo

Unfortunately if $Y$ is a sequence space the result is not
quite so clean.
We first state the corresponding theorem and then a more
general technical result which includes the theorem.

\proclaim{Theorem 3.5}  Let $Y$  be a  Banach sequence space
and  suppose  that  $(f_n)$  is a complemented unconditional
basic sequence in $Y.$
Suppose $(f_n)$ is lattice anti-Euclidean.
Then  $(f_n)$  is  equivalent  to  a  complemented  disjoint
sequence $(f'_n)$ in $Y^N$ for some natural number $N.$
\endproclaim

\proclaim{Theorem   3.6}Let   $\Cal   X$   be   a    lattice
anti-Euclidean family of Banach sequence spaces.
Then given $M>1$ there is a constant $C=C(M)$ and a  natural
number $N=N(M)$ so that the following property holds.

Suppose $Y$ is a Banach sequence space and $X\in\Cal X$ with
$\dim X=d\le \infty.$
Suppose $S\colon X\to  Y$ and $T\colon  Y\to X$ are  bounded
operators with $\|S\|,\|T\|\le M.$
Let $f_n=Se_n$  and $g_n=T^*e_n$  for $n\le  d$ and  suppose
$E_n$   are   disjoint   subsets   of   $\bold  N$  so  that
$\|f_ng_n\chi_{E_n}\|_1\ge M^{-1}.$
Then we can find subsets  $(F_{kn})$ of $\bold N$ for  $1\le
k\le N$ and $1\le n\le  d$ so that (1) $F_{kn}\subset  E_n,$
(2) for  each fixed  $k,$ the  sets $(F_{kn})$  are pairwise
disjoint, (3) for  each fixed $n,$  the sets $(F_{kn})$  are
pairwise  disjoint  and  (4)  the  disjoint  sequence $(f'_n
)_{n=1}^d$ defined  by $f_n'=(f_n\chi_{F_{kn}})_{k=1}^N$  in
$Y^N$  is  $C$-complemented  and  $C$-equivalent to the unit
vectors $(e_n)_{n=1}^d$ in $X.$ \endproclaim

\demo{Remark}Of course  if we  take $\Cal  X$ as  having one
member and $E_n=\bold N$, this implies Theorem 3.5.
However,   the   quantitative   version   will  be  of  some
importance.  \enddemo

\demo{Proof}Let $\delta=\delta(M,\Cal  X)$ be  determined as
in Lemma 3.3.
We     will     show     that     $N=[2M\delta^{-1}]$    and
$C=2K_G^2\delta^{-1}M^2N^2$ have the property claimed.

Let $h_n=|f_ng_n|\chi_{E_n}.$
Then  by  Lemma  3.3  we  can  define  an  operator $R\colon
\ell_1\to   \ell_1(c_0)$   by  $R(\zeta)=(\zeta(k)h_k)_{k=1}
^{\infty}.$
(For  notational  convenience  we  will  assume that $d=\dim
X=\infty$; minor modifications can be made if $d<\infty.$)
Now $\|R\|\le  M^2$ and  $\|R(\zeta)\|\ge \delta\|\zeta\|_1$
for all $\zeta\in\ell_1.$
We therefore  can apply  the Hahn-Banach  theorem to  find a
linear                                            functional
$\Phi=(\phi_n)_{n=1}^{\infty}\in\ell_{\infty}(\ell_1)$    so
that             $\|\Phi\|\le             1$             and
$\Phi(R\zeta)=\delta\sum_{n=1}^{\infty}\zeta(n)$   for   all
$\zeta\in\ell_1.$
In other words,
$$
\sup_k\sum_{n=1}^{\infty}|\phi_n(k)| \le 1
$$
and
$$
\sum_{k=1}^{\infty}\phi_n(k)h_n(k) =\delta
$$
for each $n.$

Let $A_n=\{k\colon |\phi_n(k)|\ge 1/(MN)\}.$
Then
$$
\sum_{k\in A_n}|\phi_n(k)|h_n(k) \ge \frac12\delta.
$$
Now          $\sum_{n=1}^{\infty}\chi_{A_n}(k)           \le
N\sum_{n=1}^{\infty} |\phi_n(k)| \le N.$
It follows that we can decompose $A_n=\cup_{k=1}^NF_{kn}$ as
a  disjoint  union  where  for  each  $1\le k\le N$ the sets
$(F_{kn})_{n=1}^{\infty}$ are disjoint.

Let $f'_n= (|f_n|\chi_{F_{kn}})_{k=1}^{N} \in Y^N$ (which we
consider as an $\ell_2$-sum).
Similarly,   let   $g'_n    =(|g_n|\chi_{F_{kn}})_{k=1}^N\in
(Y^*)^N.$ Then
$$
\langle   f'_n,g'_n\rangle   =   \sum_{k\in    A_n}h_n(k)\ge
\sum_{k\in A_n}|\phi_n(k)|h_n(k)\ge \frac12\delta.
$$

Let $\beta_n=\langle  f'_n,g'_n\rangle$ and  define $U\colon
X\to     Y^N$     and     $V\colon     Y^N\to     X$      by
$U(\xi)=\sum_{j=1}^{\infty}\xi(j)f'_j$     and      $V(\bold
y)(j)=\beta_j^{-1}\langle \bold y,g'_j\rangle$ where  $\bold
y=(y_1,\ldots,y_N).$

Suppose $\xi\in c_{00}.$
Then
$$
\|U\xi\|_{Y^N}  \le   N\|\max_{j\ge  1}|\xi(j)f_j|\|_Y   \le
N\|(\sum_{j=1}^{\infty}\xi(j)^2f_j^2)^{1/2}     \|_Y     \le
K_GMN\|\xi\|_X.
$$
From this it  quickly follows that  $U$ is well-defined  and
bounded with $\|U\|\le K_GMN.$
On the other hand if $\xi\in c_{00}$ then
$$
\langle   V\bold   y,\xi\rangle\le   2\delta^{-1}    \langle
\sum_{k=1}^N|y_k|,  \max_{j\ge  1}|\xi(j)g_j|  \rangle   \le
2K_G\delta^{-1}MN \|\bold y\|_{Y^N}
$$
so that $\|V\| \le 2K_G\delta^{-1}MN.$
Since $VU=I_X$ the proof is complete.\qed\enddemo

\demo{Remarks}It is not possible  to improve Theorem 3.5  by
replacing $Y^N$ by $Y$.
We sketch an example.
Gowers \cite{12} (cf  \cite{13}) has constructed  a sequence
space with  the property  that every  bounded operator  is a
strictly singular perturbation of a diagonal operator.
Let    $1<p<2$    and    consider    the    space    $\tilde
G=G(\ell_p^{2^n})$ (i.e. the direct sum in the sense of  $G$
of spaces $\ell_p^{2^n}.$)
The obvious basis is anti-lattice Euclidean (in fact $\tilde
G$ is $p$-concave).
However $\tilde G$ has another unconditional basis which  is
formed by taking the Haar basis in each co-ordinate.
It may be shown that the original basis is not equivalent to
a block basis of this basis.
We remark, however, that, in this example $N=2$ suffices and
we know of no example where $N=2$ does not suffice.
A  somewhat  similar  problem  is  considered by Wojtaszczyk
\cite{26} for certain types of bases in quasi-Banach spaces.

We also remark  that a continuous  analogue of Theorems  3.4
and  3.5  is  proved  by  somewhat  similar  techniques   in
\cite{14}, Theorem 8.1.
This result  which follows  from work  in \cite{16}  was, in
fact, the basis for the proof given here.
We  have  opted  however  for  a  completely  self-contained
approach.  \enddemo

\vskip1truein
\heading
4. Applications to finite-dimensional spaces
\endheading

Before stating  our first  application, let  us recall  some
definitions from \cite{7}.
Let  $\Cal  X$  be  a  family  of  finite-dimensional Banach
sequence spaces.
Suppose first each $X\in\Cal X$ is a symmetric space.
Then  we  say  the  members  of  $\Cal X$ have a {\bf unique
symmetric  basis}  if   there  is  a   function  $\psi\colon
[1,\infty)\to[1,\infty)$ so  that if  $(u_i)_{i=1}^{\dim X}$
is a normalized  $K$-symmetric basis of  some $X\in \Cal  F$
then  $(u_i)_{i=1}^{\dim  X}$  is  $\psi(K)-$equivalent   to
$(e_i)_{i=1}^{\dim X}.$

Now  consider  the  case  when  each  $X$ is not necessarily
symmetric.

Then we  say the  members of  $\Cal X$  have an  {\bf almost
(somewhat)  unique  unconditional  basis}  if  there  is   a
function $\phi\colon [1,\infty)\times  (0,1)$ so that  given
$K\ge  1$,  then  for  any  $0<\alpha<1$  (resp.  for   some
$0<\alpha=\alpha(K)<1)$ it is  true that whenever  $X\in\Cal
X$    has     a    normalized     $K$-unconditional    basis
$(u_i)_{i=1}^{\dim X}$  then there  is a  subset $\sigma$ of
$[\dim X]$ with $|\sigma|\ge \alpha\dim X$ and a one-one map
$\pi\colon \sigma\to [\dim X]$ so that  $(e_i)_{i\in\sigma}$
is              $\phi(K,\alpha)-$equivalent               to
$(u_{\pi(i)})_{i\in\sigma}.$

The  following   theorem  shows   that  any   collection  of
finite-dimensional    spaces    which    form    a   lattice
anti-Euclidean   family   (i.e.   do   not   have  uniformly
complemented  $\ell_2^n-$sublattices)  have  almost   unique
unconditional bases.
In  particular  in  any  such  class  the symmetric basis is
unique; both these results are new.
There are,  however, numerous  results of  this type  in the
literature.
It was shown by Gowers \cite{11} that the symmetric basis is
not  unique  for  the  class  of  all  symmetric spaces, but
positive results for various classes are given in  \cite{4},
\cite{7}, \cite{15} and \cite{24}.
The   problem   of   almost   or   somewhat  uniqueness  for
unconditional bases  for various  classes was  considered in
\cite{24} and \cite{7}.

\proclaim{Theorem   4.1}Let   $\Cal   X$   be   a    lattice
anti-Euclidean family of finite-dimensional sequence spaces.
Then   the   members   of   $\Cal   X$  have  almost  unique
unconditional bases.\endproclaim

In order to prove this we will need a lemma, due essentially
to Wojtaszczyk \cite{26}.
Our  statement  is  a  modification  and  we  will avoid the
language of bipartite graph theory.
Suppose $n\in\bold N$ and let $G$ be a subset of  $[n]\times
[n].$
For $i\in [n]$ let $V_i$ be the set of $j$ so that for  some
$k\in[n],$ $(i,k)$ and $(j,k)\in G.$
For $A\subset  [n]$ let  $V(A)=V^1(A)=\cup_{i\in A}V_i$  and
then   define   inductively   $V^r(A)=V(V^{r-1}(A));$    let
$V^r_i=V^r(\{i\}).$
Finally let  $G^r$ be  the set  of $(i,j)$  so that for some
$k\in [n]$ we have $(k,j)\in G$ and $k\in V^r_i.$

\proclaim{Lemma 4.2}Assume $G$ is as above.
Suppose  $(w_{ij})_{i,j\in[n]}$  are  such that:\newline (1)
$\sum_{j=1}^nw_{ij}=1$   for    $i\in   [n]$,\newline    (2)
$\sum_{i=1}^nw_{ij}=1$ for $j\in [n].$\newline
Let  $M=  \max_{1\le  j\le  n}\sum_{i=1}^n|w_{ij}|$  and let
$b=\sum_{(i,j)\notin G}|w_{ij}|.$
Then for any  $r$ there is  a subset $\sigma$  of $[n]$ with
$|\sigma|\ge  n-3b-Mr^{-1}n$  and  a  one-one map $\pi\colon
\sigma\to   [n]$   with   $(i,\pi(i))\in   G^r$   for  $i\in
\sigma.$\endproclaim

\demo{Proof}Note first that either  $V_i$ is empty or  $i\in
V_i$; let $E$ be the set of $i$ such that $V_i$ is empty.
Then $|E|=\sum_{i\in E}\sum_{j=1}^nw_{ij}\le b.$
For  any  $A$  we  have  $A\subset V(A)\cup E$ and $V(A)\cap
E=\emptyset.$
Thus the  sequence $(V^s(A))_{s=1}^n$  is increasing  and so
for every $1\le s\le r$ we have
$$
|V^s(A)|\ge |A|-b.
$$
For future reference we  let $V^0(A)=A\setminus E$ and  have
the same inequality.
Now for any $A\subset  [n]$ let $A^*=\{j\colon \exists  i\in
A,\  (i,j)\in  G\}$  and  $A^+=\{j\colon  \exists  i\in  A:\
(i,j)\in G^r\}.$
Then $A^+=(V^r(A))^*.$

Assume for some $A$ we have $|A^+|\le |A|.$
Then   there   exists   some   $0\le   s\le   r-1$  so  that
$|V^{s+1}(A)^*|\le |V^s(A)^*| +r^{-1}|A|.$
(Here we recall $V^0(A)=A\setminus E.$).
Notice $V^{s+1}(A)\supset V^s(A).$
We now compute
$$
\align               |V^{s+1}(A)|               &=\sum_{i\in
V^{s+1}(A)}\sum_{j=1}^nw_{ij}\\      &\le       b+\sum_{i\in
V^{s+1}(A)}\sum_{j\in  V^{s+1}(A)^*}w_{ij}   \\  &\le   b  +
M(|V^{s+1}(A)^*|-|V^s(A)^*|)          +           \sum_{i\in
V^{s+1}(A)}\sum_{j\in V^s(A)^*}w_{ij}\\ &\le 2b+  Mr^{-1}|A|
+ \sum_{i=1}^n\sum_{j\in V^s(A)^*}w_{ij}\\ &= 2b +Mr^{-1}|A|
+ |V^s(A)^*|.  \endalign
$$
Hence
$$
|A| \le 3b + Mr^{-1}|A| + |A^+|.
$$

Combining we have $|A^+| \ge |A|-3b-Mr^{-1}n.$
The  result   then  follows   from  Hall's   Marriage  Lemma
\cite{2}.\qed\enddemo

\demo{Proof of Theorem 4.1} It will suffice to show that for
any  $\epsilon>0$   and  $M\ge   1$  there   is  a  constant
$C=C(\epsilon,M,\Cal  X)$  with  the  property that whenever
$X,Y$ are finite-dimensional sequence spaces with  $X\in\Cal
X$  and  if  $S\colon   X\to  Y$  is  an   isomorphism  with
$\max(\|S\|,\|S^{-1}\|)\le  M$  then   there  is  a   subset
$\sigma$ of $[\dim X]$ with $|\sigma|\ge (1-\epsilon)\dim X$
and a one-one map $\pi\colon \sigma\to [\dim X]$ so that for
any $(\alpha_i)_{i\in\sigma}$ we have
$$
C^{-1}\|\sum_{i\in\sigma}\alpha_ie_i\|_X                 \le
\|\sum_{i\in\sigma}\alpha_ie_{\pi(i)}\|_Y                \le
C\|\sum_{i\in\sigma}\alpha_i e_i\|_X.
$$

Let $n=\dim X$ and let $T=S^{-1}.$
As before,  let $Se_i=f_i$  and $T^*e_i=g_i$  for $1\le i\le
n.$
Let $w_{ij}=f_i(j)g_i(j).$
Then    $\sum_{i=1}^nw_{ij}=1$     for    all     $j$    and
$\sum_{j=1}^nw_{ij}=1$ for all $i$.
Furthermore, $\sum_{i=1}^n|w_{ij}|=\||S^*e_j||Te_j|\|_1  \le
\|S^*e_j\|_{X^*}\|Te_j\|_X \le M^2.$

Now by Lemma 3.3 there exists a $\delta=\delta(\epsilon, M)$
so that if for some  subset $\tau$ of $[n]$ we  choose $0\le
h_k\le  |f_k||g_k|$  for  $k\in\tau$  so  that   $\|h_k\|\ge
\epsilon/8$ then
$$
\|\max_{k\in\tau}|\alpha_kh_k|\|_1                       \ge
\delta\sum_{k\in\tau}|\alpha_k|
$$
for all $(\alpha_k)_{k\in\tau}.$
Next choose  $\eta=\frac18\delta\epsilon M^{-2}.$

We let $G$ be the set of pairs $(i,j)$ so that  $|w_{ij}|\ge
\eta.$
Then
$$
\sum_{(i,j)\notin          G}|w_{ij}|          =\sum_{i=1}^n
\|f_ig_i\chi_{A_i}\|_1
$$
where $A_i=\{j\colon |w_{ij}|<\eta.$
Let   $h_i=|f_ig_i|\chi_{A_i}$   and   let   $\tau=\{i\colon
\|h_i\|_1\ge \frac18\epsilon\}.$
Then
$$
\delta|\tau|\le \|\max_{i\in\tau}h_i\|_1 \le n\eta
$$
so that $|\tau|\le \frac18\epsilon M^{-2}n.$
Hence $\sum_{i\in\tau}\|h_i\|_1 \le \frac18\epsilon n.$
However $\sum_{i\notin\tau}\|h_i\|_1 \le \frac18\epsilon n.$
Combining we have
$$
\sum_{(i,j)\notin G}|w_{ij}|\le \frac14\epsilon n.
$$

We can now apply Lemma 4.2.
We  choose  $r=[4M^2\epsilon]+1$  so  that  we have a subset
$\sigma$  of  $[n]$  with  $|\sigma|\ge (1-\epsilon)n$ and a
one-one map $\pi\colon \sigma\to [n]$ so that $(i,\pi(i))\in
G^r$ for $i\in\sigma.$

Note  that  if   $(i,j)\in  G$  then   $|f_i(j)|,|g_i(j)|\ge
\delta/M.$
Now by Krivine's theorem if $(\alpha_i)_{i=1}^n$ are scalars
then
$$
\|(\sum_{i=1}^n\alpha_i^2f_i^2)^{1/2}\|_Y                \le
K_GM\|\sum_{i=1}^n\alpha_ie_i\|_X.
$$
Hence
$$
\|(\sum_{(i,j)\in       G}\alpha_i^2e_j)^{1/2}\|_Y       \le
K_GM^2\delta^{-1}\|\sum_{i=1}^n\alpha_ie_i\|_X.
$$
Similarly, a dual argument gives that
$$
\|(\sum_{(i,j)\in       G}\alpha_j^2e_i)^{1/2}\|_X       \le
K_GM^2\delta^{-1}\| \sum_{j=1}^n\alpha_je_j\|_Y.
$$
Iterating     these     conditions     gives     that     if
$C=(K_G^2M^2\delta^{-1})^{r+1}$ then
$$
\|(\sum_{(i,j)\in      G^r}\alpha_i^2e_j)^{1/2}\|_Y      \le
C\|\sum_{i=1}^n\alpha_ie_i\|_X.
$$
It follows that
$$
\|\sum_{i\in\sigma}\alpha_ie_{\pi(i)}\|_Y                \le
C\|\sum_{i\in\sigma}\alpha_ie_i\|_X.
$$
Similarly we have
$$
\|(\sum_{(i,j)\in      G^r}\alpha_j^2e_i)^{1/2}\|_Y      \le
C\|\sum_{j=1}^n\alpha_je_j\|_X,
$$
so that
$$
\|\sum_{i\in\sigma}\alpha_ie_i\|_X                       \le
C\|\sum_{i\in\sigma}\alpha_ie_{\pi(i)}\|_Y
$$
so that the result follows.\qed\enddemo

\vskip1truein

\heading
5. Right- and left-dominant spaces
\endheading

Let $X$ be a sequence space.
We will  say that  $X$ is  {\bf left-dominant  with constant
$\gamma\ge  1$}   if  whenever   $(u_1,u_2,\ldots,u_n)$  and
$(v_1,\ldots,v_n)$ are  two disjoint  sequences in  $c_{00}$
with  $\|u_k\|_X\ge  \|v_k\|_X$  and  such  that $\supp v_k>
\supp u_k$ for $1\le k\le n$ then $\|\sum_{k=1}^nv_k\|_X \le
\gamma\|\sum_{k=1}^nu_k\|_X.$
Similarly, we will say that $X$ is {\bf right-dominant  with
constant  $\gamma$}   if  whenever   $(u_1,\ldots,u_n)$  and
$(v_1,\ldots,v_n)$   are   two   disjoint   sequences   with
$\|u_k\|_X\le  \|v_k\|_X$  and  $\supp  v_k>  \supp u_k$ for
$1\le     k\le     n$     then     $\|\sum_{k=1}^nu_k\|_X\le
\gamma\|\sum_{k=1}^nv_k\|_X.$

We will refer to any normalized unconditional basic sequence
as being left- or right-dominant according as the associated
sequence space is left- or right-dominant.

\proclaim{Lemma  5.1}$X$  is  left-dominant  with   constant
$\gamma$  if  and  only  if  $X^*$  is  right-dominant  with
constant $\gamma.$\endproclaim

\demo{Proof}Let us prove that  if $X$ is left-dominant  then
$X^*$ is right-dominant; the other direction is similar.
Suppose  $(u_1^*,\ldots,u_n^*)$  and  $(v_1^*,\ldots,v_n^*)$
are  two  sequences  in  $c_{00}$  with  $\|u_k^*\|_{X^*}\le
\|v_k^*\|_{X^*}$ for  $1\le k\le  n$ and  $\supp u_k^*<\supp
v_k^*.$
There exists $x\in X$ supported on $\cup_{i=1}^n\supp u_i^*$
with            $\|x\|=1$            and            $\langle
x,\sum_{i=1}^nu_i^*\rangle=\|\sum_{i=1}^nu_i^*\|_{X^*}.$
Let  $x=\sum_{i=1}^nu_i$   where  $\supp   u_i\subset  \supp
u_i^*.$
Next pick $v_i$ of norm one with support contained in $\supp
v^*_i$ so that $\langle v_i,v_i^*\rangle=\|v_i^*\|_{X^*}.$
Finally let $y=\sum_{i=1}^n\|u_i\|_Xv_i.$
Then $\|y\|\le \gamma \|x\|=\gamma.$
Also
$$
\align          \langle           y,\sum_{i=1}^nv_i^*\rangle
&=\sum_{i=1}^n\|u_i\|_X\|v_i^*\|_{X^*}\\                &\ge
\sum_{i=1}^n\langle      u_i,u_i^*\rangle      \\       &\ge
\|\sum_{i=1}^nu_i^*\|_{X^*}.\qed\endalign
$$
\enddemo

If $N$  is a  natural number  we denote  by $X_N$  the space
$X[N+1,\infty)$  of  all  $x\in  X$  such that $x(k)=0$ when
$k\le N.$

\proclaim{Lemma   5.2}   Suppose   $X$   is   a  left-(resp.
right-)dominant sequence space.
Suppose  $1\le  p\le  \infty$  and  $\ell_p$  is  disjointly
finitely representable in $X$.
Then  :\newline  (1)  $X$  satisfies  a lower-(resp. upper-)
$p$-estimate .
\newline  (2)  There  is  a  constant  $K$  so  that for any
$n\in\bold  N,$  there  exists  $N\in\bold  N$ so that $X_N$
satisfies an  upper- (resp.  lower-) estimate  with constant
$K$ on $n$ vectors.  \endproclaim

\demo{Proof}We consider only the case of a left-dominant
space, and assume that $\ell_p$ is $C$-disjointly
representable in $X$ (actually by Krivine's theorem
\cite{19} we could suppose $C=1$).
For notational convenience suppose $p<\infty.$
Suppose $x_1,\ldots,x_n$ are disjoint in $c_{00}.$
Then there exist $y_1,\ldots,y_n$ disjoint in $X$ with
$\max_k\supp x_k<\supp y_j$ for $1\le j\le n$ such that
$\|y_j\|=\|x_j\|$ and $\|\sum_{j=1}^n y_j\| \ge (2C)^{-1}
(\sum_{j=1}^n\|x_j\|^p)^{1/p}$.
Thus $\|\sum_{j=1}^nx_j\|\ge (2C\gamma)^{-1}
(\sum_{j=1}^n\|x_j\|^p)^{1/p}$.
Conversely, if we fix $n$ and choose any $y_1,\ldots,y_n$
normalized, disjoint and $2C$-equivalent to an
$\ell_p^n$-basis then if $N=\max_k\supp y_k$ and
$x_1,\ldots,x_n$ are disjoint in $X_N$ then $
\|\sum_{j=1}^nx_j\| \le \gamma \|\sum_{j=1}^n\|x_j\|y_j\|
\le (2C\gamma)\|\sum_{j=1}^n\|x_j\|^p)^{1/p}.$ \qed\enddemo

It follows from Lemma 5.2 and Krivine's theorem that if $X$
is left- or right-dominant then there is exactly one
$r=r(X)$ so that $\ell_r$ is disjointly finitely
representable in $X$.
Let us call $r$ the {\bf index} of $X.$
If $X$ is right-dominant and of index $\infty$ then clearly
$X=c_0$ while if $X$ is left-dominant of index $1$ then
$X=\ell_1.$
A right-dominant space of finite index has a nontrivial
lower estimate and so can realized as the dual of
left-dominant space of index greater than one.

Notice that it also follows from Lemma 5.2 that every left-
or right-dominant sequence space is an asymptotically
$\ell_r$space where $r=r(X)$ (cf.  \cite{23}, p. 221).
We now turn our attention to the problem of deciding when a
left- or right-dominant space is sufficiently Euclidean.

\proclaim{Proposition 5.3}Let $X$ be a left- or
right-dominant sequence space.  Then $X$ is sufficiently
Euclidean if and only if $1<r(X)<\infty.$ \endproclaim

\demo{Proof}Let $\Cal U$ be a nonprincipal ultrafilter on
the natural numbers and let $X_n=X[n,\infty).$
Let $Y$ be the ultraproduct $\ell_{\infty}(X_n)/c_{0,\Cal
U}(X_n)$ where $c_{0,\Cal U}(X_n)$ consists of all sequences
$(x_n)\in \ell_{\infty}(X_n)$ with $\lim_{n\in\Cal
U}\|x_n\|=0.$
Then $X$ is sufficiently Euclidean if and only if $\ell_2$
embeds complementably into $Y.$
Assume $X$ is left- or right-dominant with index $r$.
Then $Y$ is a Banach lattice with an upper and lower
$r$-estimate.  This implies $Y$ is isomorphic to an abstract
$L_r$-space and so the result follows.\qed\enddemo

\proclaim{Proposition 5.4}Suppose $(u_n)$ is a left- (resp.
right-) dominant basis and that $\pi$ is a permutation of
the natural numbers such that $(u_{\pi(n)})$ is also left-
(resp. right-)dominant.
Then there is a constant $C$ such that for any $\alpha\in
c_{00}$
$$
\|\sum_{k=1}^{\infty}\alpha_ku_{2k}\| \le C
\|\sum_{k=1}^{\infty}\alpha_ku_{\pi(k)}\|
$$
(respectively,
$$
\|\sum_{k=1}^{\infty}\alpha_ku_{2k}\| \ge C^{-1}
\|\sum_{k=1}^{\infty}\alpha_ku_{\pi(k)}\|.)
$$
\endproclaim

\demo{Proof}We treat only the left-dominant case.
Define a sequence $(s_n)$ inductively as follows.
Let $s_1=1$ and then let $s_n$ be the least $m$ so that
$m\in \pi\{n,n+1,\ldots\}\setminus\{s_1,\ldots,s_{n-1}\}.$
Note that $s_n\le 2n-1<2n$ and that $(s_n)$ increases.
Further $s_n=\pi(r_n)$ where $r_n\ge n$.
Hence $(u_{\pi(n)})$ dominates $(u_{s_n})$ which in turn
dominates $(u_{2n}).$
This establishes the result.\qed\enddemo

\proclaim{Proposition 5.5}Let $(u_n)$ be a left- (resp.
right-)dominant basis of a Banach space $X.$
In order that $(u_n)$ be equivalent to its square it is
necessary and sufficient that $(u_n)$ be equivalent to
$(u_{2n}).$ \endproclaim

\demo{Remark}Clearly $(u_n)$ is equivalent to $(u_{2n})$ if
and only if $(u_n)$ is equivalent to $(u_{Nn})$ for any $N$
in view of the dominance assumption.\enddemo

\demo{Proof}Consider the left-dominant case and assume
$(u_n)$ is equivalent to its square.
Let $(v_n)$ be the natural basis of $X^4$ with $(v_{4n-j})$
equivalent to $(u_n),$ for $0\le j\le 3.$ Then $(v_n)$ is
left-dominant.
Since some permutation of $(v_n)$ is equivalent to $(u_n)$
we have that $(u_n)$ dominates $(v_{2n})$ and hence
$(u_{2n})$ dominates $(u_n)$.
This implies $(u_n)$ and $(u_{2n})$ are equivalent.
The other case is similar and the other direction is
trivial.\qed\enddemo

\proclaim{Theorem 5.6}Let $X$ be a left- or right-dominant
separable sequence space and that $(u_n)$ is a complemented
normalized disjoint sequence in $X.$
Then $(u_n)$ is permutatively equivalent to a subsequence of
the canonical unconditional basis $(e_n)$ of $X.$
\endproclaim

\demo{Proof}Let us assume the basis is left-dominant; the
case of a right-dominant basis is almost identical.
We can assume the dual functionals $(u_n^*)$ in $X^*$ have
the same support as $(u_n).$ Let
$f_n=|u_n||u_n^*|\in\ell_1.$
For each $n$ pick $k_n\in\supp u_n$ so that
$\|f_ne_{[1,k_n]}\|_1\ge \frac12$ and
$\|f_ne_{[k_n,\infty)}\|_1\ge \frac12.$
Let $v_n=u_ne_{[1,k_n]}$ and $w_n=u_ne_{[k_n,\infty)}.$
Now we argue that $(v_n)$ and $(w_n)$ are both equivalent to
$(u_n).$
Indeed the operator $Tx=\sum_{n=1}^{\infty}\langle
x,|u_n^*|\rangle u_n$ is easily seen to be bounded on $X$.
We have $T|v_n|=\alpha_n u_n,$ and $T|w_n|=\beta_nu_n$ where
$\alpha_n,\beta_n\ge 1/2.$
It follows that both $(|v_n|)$ and $(|w_n|)$ are equivalent
to $(u_n)$ and the desired conclusion follows.

Now if $X$ is left-dominant $(v_n)$ dominates $(e_{k_n})$;
to see this just note that $(v_ne_{[1,k_n)})$ dominates
$(\|v_ne_{[1,k_n)}\|e_{k_n})$.  Similarly, $(e_{k_n})$
dominates $(w_n)$.
Thus $(u_n)$ is equivalent to $(e_{k_n}).$ \qed \enddemo

\proclaim{Theorem 5.7}Let $X$ be a separable left- or right-
dominant sequence space.
Suppose that $r(X)=1$ or $r(X)=\infty$ and that $(e_n)$ and
$(e_{2n})$ are equivalent.
Then every complemented normalized unconditional basic
sequence is equivalent to a subsequence of the basis and $X$
has a unique unconditional basis.\endproclaim

\demo{Proof}We assume $X$ left-dominant.
Let $(u_n)$ be any complemented normalized unconditional
basic sequence in $X$.
Then $(u_n)$ is anti-lattice Euclidean by Proposition 5.3
and so by Theorem 3.5 and the hypothesis on $X$, $(u_n)$ is
equivalent to a complemented disjoint sequence in $X.$
By Theorem 5.6, this implies that $(u_n)$ is permutatively
equivalent to a subsequence $(e_{k_n})$ of $(e_n)$.

We now restrict to the case when $(u_n)$ is an unconditional
basis of $X.$
Applying Theorem 3.5 again we see that $(e_n)$ is equivalent
to a complemented disjoint sequence in the $N$-fold sum
basis $(u_n)^N$ of $X^N$.
Now $(e_{k_n})^N$ arranged in the obvious order is also a
left-dominant basis.
Here the ``obvious order'' $(f_n)$ is to take $f_{N(j-1)+s}$
to be $(0,\ldots,0,e_{k_j},0,\ldots)\in X^N$, where
$e_{k_j}$ is in the $s$th. co-ordinate.
Hence $(e_n)$ is permutatively equivalent to a subset of
$(e_{k_n})^N.$
However $(e_{k_n})^N$ is permutatively equivalent to a
subset of $(e_n)^N$ which is permutatively equivalent to
$(e_n).$
By the Cantor-Bernstein principle \cite{27}, this means that
$(e_{k_n})^N$ is permutatively equivalent to $(e_n).$

Now $(f_k)$ dominates $(e_{2k})$ by Proposition 5.4 and
similarly $(e_k)$ dominates $(f_{2k}).$
Hence, since $(e_k)$ and $(e_{4k})$ are equivalent we have
that $(f_{2k})$ is equivalent to $(e_k).$
Now $(f_{2n-1})_{n\ge 1} $ is dominated by
$(f_1,f_2,f_4,\ldots)$ and dominates $(f_1,f_4,f_8,\ldots)$
and thus is also equivalent to $(f_{2n}).$
Hence $(f_n)$ is equivalent to $(e_n).$
Now $f_{Nn}$ is equivalent to $e_{Nn}$ and hence to $(e_n)$.
Thus $(e_{k_n})$ is equivalent to $(e_n).$
The result now follows.\qed\enddemo

\demo{Remarks} There is a natural question here, which is
also suggested by the work of Wojtaszczyk \cite{26}.
Suppose $(x_n)$ and $(y_n)$ are two unconditional bases
whose squares are permutatively equivalent; does it follow
that $(x_n)$ and $(y_n)$ are permutatively equivalent?
The corresponding Banach space problem has a negative
solution.
An example of Gowers \cite{13} shows that there is a Banach
space $X$ so that $X$ and $X^2$ are not isomorphic but $X^2$
and $X^4$ are isomorphic.  \enddemo

\proclaim{Theorem 5.8}Suppose that $1<p_n<\infty$ for all
$n$ and that $p_n\downarrow 1.$
Suppose that for some constant $a>0$
$$
\frac1{p_{2n}}-\frac1{p_n}\le \frac a{\log n}
$$
for $n\ge 2.$
Then the Nakano space $\ell(p_n)$ has a unique unconditional
basis.

Similarly if $p_n\uparrow \infty$ and for some constant
$a>0$
$$
\frac1{p_n}-\frac1{p_{2n}}\le \frac a{\log n}
$$
for $n\ge 2,$ then $h(p_n)$ has a unique unconditional
basis.\endproclaim

\demo{Proof}If $p_n\downarrow 1$ then $X=\ell(p_n)$ has a
right-dominant basis with $r(X)=1.$
The assumption that $\frac1{p_n}-\frac1{p_{2n}}=O((\log
n)^{-1})$ implies that the basic sequences $(e_n)$ and
$(e_{2n})$ are equivalent by an old result of Simons
\cite{25}.
The second case is similar (or dual).\qed\enddemo

\proclaim{Theorem 5.9}Let $p_n\downarrow 1$ be such that
$p_n= 1+O(\log n)^{-1}.$
Assume $(N_n)$ is an increasing sequence of natural numbers
such that $\frac{1}{p_{n+1}}-\frac1{p_n} = O((\log
N_n)^{-1})$ and $\inf_n{N_{n+1}/N_n}>1.$
Then $\ell_1(\ell_{p_n}^{N_n})$ has a unique unconditional
basis.\endproclaim

\demo{Proof}Let $q_n$ be the sequence obtained by writing
out $p_1,$ $N_1$ times, $p_2$ $N_2$ times etc.
It is clear that the Nakano space $\ell(q_n)$ can be written
as vector-valued Nakano space $\ell(p_n)(\ell_{p_n}^{N_n})$.
But by Simons's theorem \cite{25} $\ell(p_n)=\ell_1.$
Let $M_n=N_1+\cdots+N_n.$
Then $M_n\le cN_n$ for some $c$.  If $M_{n-1}< k\le M_n$ and
$M_{m-1}<2k\le M_m$ we have that either $m-n\le 1$ or
$(m-n-1)N_n\le k.$
Thus $(m-n)$ is bounded independent of $k$ and so
$\frac1{q_{2k}}-\frac1{q_{k}}=O((\log N_n)^{-1})=O(\log
k)^{-1}.$
The result follows from Theorem 5.8.\qed\enddemo

\demo{Remark} We do not know if Theorem 5.9 holds for any
space $\ell_1(\ell_{p_n}^{N_n})$ where $p_n\downarrow
1.$\enddemo

Our final example of this section is the now classical
Tsirelson space.  We refer to \cite{8} for full details of
this space.
We recall that the Tsirelson norm $\|\,\|_T$ on $c_{00}$ is
the minimal norm satisfying $\|x\|_T \ge \|x\|_{\infty}$ and
$$
\frac12\sum_{j=1}^n\|x_j\|_T \le \|\sum_{j=1}^nx_j\|_T
$$
whenever $n\le\supp x_1<\supp x_2<\cdots<\supp x_n.$
Tsirelson space is the sequence space $T$ obtained by
completing $c_{00}$ with respect to this norm.
This space is the dual of the original Tsirelson space.
We will need an alternative norm $\|\,\|^\#_T$ which is
defined to be the least norm satisfying $\|x\|_T^\#\ge
\|x\|_{\infty}$ and
$$
\frac12\sum_{j=1}^{2n}\|x_j\|_T^\# \le
\|\sum_{j=1}^{2n}x_j\|_T^\#
$$
whenever $x_1,\ldots,x_{2n}$ are disjoint and $n\le\supp
x_j$ for $j=1,2,\ldots,2n.$

\proclaim{Lemma 5.10} For $x\in c_{00}$ we have $\|x\|_T \le
\|x\|_T^{\#}\le 4\|x\|_T.$\endproclaim

\demo{Proof}Let $\|\,\|'_T$ be the least norm on $c_{00}$ so
that $\|x\|'_T \ge \|x\|_{\infty}$ and
$$
\frac12\sum_{j=1}^{2n}\|x_j\|'_T \le
\|\sum_{j=1}^{2n}x_j\|'_T
$$
whenever $n\le\supp x_1<\supp x_2<\cdots<\supp x_{2n}.$
By \cite{5} we have $\|x\|_T \le \|x\|'_T \le 2\|x\|_T$ and
by \cite{1} we have $\|x\|'_T \le \|x\|_T^{\#}\le
2\|x\|'_T.$\qed\enddemo

We now prove that Tsirelson space is right-dominant.
This result was stated without proof in \cite{3}, and
generalizes Lemma II.1 of \cite{8}.

\proclaim{Lemma 5.11}Suppose $x_1,\ldots, x_N$ are disjoint
in $c_{00}$ and let $a_k=\max\supp x_k.$
Then
$$
\|\sum_{k=1}^Nx_k\|_T \le 4 \|\sum_{k=1}^N
\|x_k\|_Te_{a_k}\|_T^{\#}.
$$
\endproclaim

\demo{Proof} Indeed if this inequality is false there exist
disjoint $x_1,\ldots,x_N$ with support
$\supp(x_1+\cdots+x_N)$ of minimal cardinality such that
$$
\|\sum_{k=1}^Nx_k\|_T >
4\|\sum_{k=1}^N\|x_k\|_Te_{a_k}\|_T^{\#}.
$$
Let $x=\sum_{k=1}^Nx_k.$
Then clearly $\|x\|_T> \|x\|_{\infty}.$
Hence there exists $n\ge 2,$ and finite intervals $n\le
E_1<E_2<\cdots<E_n$ so that $E_kx\neq 0$ for
$k=1,2,\ldots,n$ and
$$
\|x\|_T =\frac12\sum_{k=1}^n \|E_kx\|_T.
$$
Using the minimal cardinality of $x$ we have that
$\cup_{k=1}^nE_k$ contains $\supp x.$
Note first that for any $j$ we have
$$
\|x_j\|_T \ge \frac12\sum_{k=1}^n\|E_kx_j\|_T.
$$

Now let $G_k=\{j\colon \supp x_j\subset E_k\}$ and let $H_k$
be the set of $j$ so that $a_j\in E_k$ but $j\notin G_k.$

Then $|H_k| < \min E_k$ so that
$$
 \|\sum_{j\in H_k}\|x_j\|_Te_{a_j}\|_T \ge \frac12 \sum_{j\in
H_k}\|x_j\|_T.
$$
Thus
$$
 \frac12\sum_{j\in H_k}\sum_{l=1}^n\|E_lx_j\|_T  \le 2\|\sum_{j\in
H_k}\|x_j\|_Te_{a_j}\|_T.
$$

Also, by our minimality assumption we have
$$
\frac12 \|E_k\sum_{j\in G_k}x_j\|_T \le 2\|\sum_{j\in
G_k}\|x_j\|_Te_{a_j}\|_T^{\#}.
$$

Combining these statements we obtain
$$
\align \frac12\sum_{k=1}^n\|E_kx\|_T &\le
2\sum_{k=1}^n\|\sum_{j\in G_k}\|x_j\|_Te_{a_j}\|_T^{\#}+
2\sum_{k=1}^n\|\sum_{j\in H_k}\|x_j\|_Te_{a_j}\|_T^{\#}\\
&\le 4\|\sum_{j=1}^N\|x_j\|_Te_{a_j}\|_T^{\#} \endalign
$$
as required.\qed\enddemo

\proclaim{Proposition 5.12}Tsirelson space $T$ is
right-dominant.  \endproclaim

\demo{Proof}Suppose that $(x_j)_{j=1}^n,(y_j)_{j=1}^n$ are
two disjoint sequences with $x_j,y_j\in c_{00}$ and
$\|x_j\|_T\le \|y_j\|_T$ and $\supp x_j<\supp y_j$ for $1\le
j\le n.$
Let $a_j=\max \supp x_j.$
Then
$$
\|\sum_{j=1}^nx_j\|_T \le
4\|\sum_{j=1}^n\|x_j\|_Te_{a_j}\|_T^{\#}
$$
The proof of Lemma II.1 of \cite{8} works for the norm
$\|\,\|_T^{\#}$ with only notational changes and yields that
$$
\|\sum_{j=1}^n\|x_j\|_Te_{a_j}\|_T^{\#} \le
\|\sum_{j=1}^ny_j\|_T^{\#}
$$
since $\|x_j\|_T\le \|y_j\|_T^{\#}.$
Combining we obtain that
$$
\|\sum_{j=1}^nx_j\|_T\le 16\|\sum_{j=1}^ny_j\|_T.
$$
Thus $T$ is right-dominant.\qed\enddemo

\demo{Remark} Of course this implies that $p-$convexified
Tsirelson $T^{(p)}$ also is right-dominant for
$1<p<\infty.$\enddemo

\proclaim{Theorem 5.13}Tsirelson space and its dual have
unique unconditional bases.\endproclaim

\demo{Proof}We have $T$ right-dominant and clearly $r(T)=1.$
We need only observe that the canonical basis $(e_n)$ is
equivalent to $(e_{2n})$ in $T$ (\cite{8} p. 14) and apply
Theorem 5.7.  \qed\enddemo

Theorem 5.13 answers a question in \cite{3}, where it is
shown that convexified Tsirelson $T^{(2)}$ has a unique
unconditional basis.
In fact much more is true as with $T^{(2)}$ (cf.  Theorem
7.9 of \cite{3}).
In fact one could prove Theorems 5.13 and 5.14 directly from
Theorem 3.5, by using known results, but it seems more
natural to invoke the theory of right-dominant bases as
here.

\proclaim{Theorem 5.14} Every complemented subspace of $T$
with an unconditional basis has a unique unconditional
basis.\endproclaim

\demo{Proof}By Theorem 5.7, every complemented normalized
unconditional basic sequence is equivalent to a subsequence
of the canonical basis.
The result follows in the same way as the preceding result,
since every subsequence of the basis is right-dominant and
equivalent to its square (\cite{8} p. 14).\qed\enddemo

\vskip1truein

\heading
6. Further examples: Orlicz sequence spaces
\endheading

In this section we construct some examples of spaces with
unique unconditional basis but such that some complemented
subspace fails to have unique unconditional basis.

Let $F$ be an Orlicz function satisfying the
$\Delta_2-$condition, normalized such that $F(1)=1.$
If we set $\phi(\tau)=sF'(s)/F(s)$ where $s=e^{-\tau}$ then
we can write $F$ in the form
$$
F(t) =\exp\left(-\int_0^{\log(1/t)}\phi(\tau)\,d\tau\right)
$$ for $0<t\le 1.$
It will be convenient to let
$\Phi(u)=\int_0^u\phi(\tau)d\tau$ for $s\ge 0.$

\proclaim{Lemma 6.1}Suppose $x_1,\ldots,x_n$ are disjoint in
$\ell_F$, and satisfy $\|x_k\|_{\ell_F}=1$.
Suppose $e^{-\tau_k}=\|x_k\|_{\infty}$ and let
$q_k=\sup_{\tau\ge \tau_k}\phi(\tau)$ and $r_k=\inf_{\tau\ge
\tau_k}\phi(\tau).$
Then for any $a_1,\ldots,a_n\in\bold R$ with
$\|\sum_{k=1}^na_kx_k\|_{\ell_F}=1$ we have
$$
\sum_{k=1}^n|a_k|^{q_k}\le 1\le \sum_{k=1}^n|a_k|^{r_k}.
$$
\endproclaim

\demo{Proof}The proof is essentially trivial.
We need only observe that if $j\in\supp x_k,$
$$
F(|a_kx_k(j)|)=
\exp(\Phi(\log|x_k(j)|^{-1})-\Phi(
\log|a_kx_k(j)|^{-1})F(|a_kx_k(j)|),
$$
and that, since $|a_k|\le 1$ and $|x_k(j)|\le e^{-\tau_k},$
$$
r_k\log |a_k|^{-1}\le \Phi(\log |a_kx_k(j)|^{-1})-\Phi(\log
|x_k(j)|^{-1})\le q_k\log |a_k|^{-1}.\qed
$$
\enddemo

\proclaim{Lemma 6.2}If $\lim_{t\to\infty}\phi(t)=1$ then the
Orlicz sequence space $\ell_F$ is
anti-Euclidean.\endproclaim

\demo{Proof}Note first that $\ell_F$ has cotype 2.
Assume that for some $M$ and every $n$ there exist operators
$S_n\colon \ell_2^{2n}\to \ell_F$ and $T_n\colon
\ell_F\to\ell_2^{2n}$ so that $T_nS_n$ is the identity on
$\ell_2^{2n}$, $\|S_n\|=1$ and $\|T_n\|\le 1.$

For fixed $n,$ we may pick by induction an orthonormal basis
$(f_k)_{k=1}^{2n}$ so that if $v\in [f_k]_{k=j+1}^{2n}$ then
$\|S_nv\|_{\infty}\le \|S_nf_j\|_{\infty}.$
Let $H_n=[f_k]_{k=n+1}^{2n}.$
Then if $v\in H_n,$ $\|S_nv\|_{\infty} \le
\|S_nf_{n+1}\|_{\infty}=\alpha_n,$ say.
For fixed $k\in \bold N$ we have $\sum_{j\in
E}|S_nf_j(k)|\le M|E|^{1/2},$ when $E\subset [n].$
It follows that if $E_k=\{j\in [n]\colon |S_nf_j(k)|\ge
\alpha_n\}$ then $|E_k| \le M^2\alpha_n^{-2}.$
But then
$$
\align \alpha_n F^{-1}(M^2\alpha_n^{-2}n^{-1})^{-1}&\le
\|(\sum_{j=1}^n |Sf_j|^2)^{1/2}\|_{\ell_F}\\ &\le \sqrt n.
\endalign
$$
Hence
$$
F(\alpha_n/\sqrt n) \le M^2\alpha_n^{-2}n^{-1}.
$$

Now if $H_n=[f_j]_{j=n+1}^{\infty}$ then
$\|S_nv\|_{\infty}\le \alpha_n$ for $v\in H_n.$
It follows from the equation above that as $n\to \infty$ we
have $\lim\alpha_n=0.$
Now let $\Cal U$ be a nontrivial ultrafilter on the natural
numbers.
Consider the ultraproduct $\ell_{\infty}(\ell_F)/c_{0,\Cal
U}(\ell_F)$ and the closed subspace thereof $Z_{\Cal
U}=Z/c_{0,\Cal U}(\ell_F)$ where $Z$ is the set of sequences
$(x_n)$ with $\lim\|x_n\|_{\infty}=0.$
Then $Z_{\Cal U}$ must contain a complemented Hilbert space.
However $Z$, as a Banach lattice, is an abstract L-space.
This follows immediately from Lemma 6.1.
Thus we have a contradiction.\qed\enddemo

The Orlicz space  $\ell_F$ has a symmetric basis and therefore every
sequence of constant coefficient blocks is a complemented unconditional
basic sequence.
Each such sequence is equivalent to the canonical basis in a
modular or Orlicz-Musielak sequence space $\ell_F[s_n]$ of
all sequences $x$ such that $
\sum_{j=1}^{\infty}F_{s_j}(|x(j)|)<\infty$ where
$F_s(t)=F(st)/F(s).$
Conversely the canonical basis of every such modular
sequence space $\ell_F[s_n]$ is equivalent a sequence of
constant coefficient blocks.
If $(s_n)$ fails to converge to $0$ then (cf \cite{19},
Proposition 3.a.5, p. 117) $\ell_F[s_n]$ is isomorphic to
$\ell_F.$
If $\lim_{n\to\infty}s_n=0$ then we can suppose that $(s_n)$
is monotonically decreasing.

Let us say that $F$ is {\bf multiplicatively convex or
m-convex} if it satisfies the condition that
$F(s^{\theta}t^{1-\theta})\le F(s)^{\theta}F(t)^{1-\theta}$
whenever $0<s,t,\theta<1.$
In this case it is clear that $\Phi$ is concave and that
$\phi$ is monotonically decreasing.

Now if $F$ is m-convex and $(s_n)$ is a monotone decreasing
sequence it is easy to see that if $\alpha\in c_{00}$ and
$(r_k)$ is an increasing sequence of natural numbers so that
$r_k\ge k$ for all $k$ then
$\|\sum_{k=1}^{\infty}\alpha_ke_{r_k}\|_{\ell_F[s_n]}\ge
\|\sum_{k=1}^{\infty}\alpha e_k\|_{\ell_F[s_n]}$.
Thus there is a weak form of dominance for the canonical
basis of $\ell_F[\tau_n].$
Based on these observations we can repeat the arguments of
Propositions 5.4 and 5.5, which only require this weakened
version, to obtain the following:

\proclaim{Lemma 6.3}Suppose $F$ is m-convex and that $(s_n)$
is sequence with $0<s_n\le 1$ and $s_n\downarrow 0.$
The canonical basis $(e_n)$ of $\ell_F[s_n]$ is equivalent
to its square if and only if $(e_n)$ is equivalent to
$(e_{2n}).$\endproclaim

\proclaim{Lemma 6.4}Suppose $F$ is m-convex.
Suppose $(s_n)_{n=1}^{\infty}$ is a monotone decreasing
sequence with $0<s_n\le 1,$ and that $(u_n)$ is a
complemented normalized disjoint sequence in $\ell_F[s_n]$.
Then there is a permutation $\pi$ of $\bold N$ and a
sequence $(s'_n)_{n=1}^{\infty}$ satisfying $0<s'_n\le s_n$
and such that $(u_{\pi(n)})$ is equivalent to the unit
vector basis of $\ell_F[s'_n].$
If in addition $\lim s_n=0$ we may suppose that $(s'_n)$ is
also decreasing.\endproclaim

\demo{Remark}If we take $s_n=1$ for all $n,$ we obtain the
fact that every complemented block basis in $\ell_F$ is
equivalent to a constant coefficient block basic
sequence.\enddemo

\demo{Proof}The proof is standard.
Suppose $(u_n^*)$ are the dual functionals and that
$f_n=|u_n||u_n^*|.$
Let $r_n=\max_{k\in \supp u_n} s_k.$
Pick $s''_n$ so that if $A_n=\{k\colon |u_n(k)|s_k\ge
s''_n\}$ and $B_n=\{k\colon |u_n(k)|s_k\le s''_n\}$ then
$\|f_ne_{A_n}\|_1,\|f_ne_{B_n}\|_1\ge \frac12.$
Then $(u_n)$ is equivalent to both the sequences
$(u_ne_{A_n})$ and $(u_ne_{B_n}).$ However, since $\phi$ is
monotone decreasing then for $k\in A_n$ and any $0\le
\lambda\le 1$ we have
$$
\frac{F(\lambda s_n |u_n(k)|)}{F(s_n|u_n(k)|)} \le
\frac{F(\lambda s''_n)}{F(s''_n)}.
$$
This implies that $(u_ne_{A_n})$ is dominated by the unit
vector basis of $\ell_F[s''_n].$
A similar argument with $B_n$ gives that $(u_n)$ is
equivalent to $\ell_F[s''_n].$

To complete the proof suppose $s> 0$ and observe that
$|\{n\colon s_n\ge s\}|\ge |\{n\colon r_n \ge s\}|\ge
|\{n\colon s''_n\le \tau\}|$ so that we can permute
$(s''_n)$ to form a sequence $(s'_n)$ with the desired
properties.  \qed\enddemo

\proclaim{Theorem 6.5}Suppose $F$ is m-convex and
$\lim_{t\to\infty}\phi(t)=1.$
Suppose $Z$ is a complemented subspace of $\ell_F$ with an
unconditional basis equivalent to its square, and such that
$Z$ is not isomorphic to $\ell_F.$
Then $Z$ has a unique unconditional basis.  \endproclaim

\demo{Proof}In fact the given unconditional basis is
equivalent to the canonical basis of $\ell_F[s_n]$ where
$s_n \downarrow 0$ and $(e_n)$ is equivalent to $(e_{2n}).$
By Lemmas 6.2 and 6.4 and Theorem 3.5 we see that any other
unconditional basis is permutatively equivalent to the unit
vector basis of $\ell_F[s'_n]$ where $s'_n\downarrow 0$ is
increasing and $s'_n\le s_n.$
But then, we can similarly find an integer $N$ so that the
original basis is permutatively equivalent to a complemented
disjoint sequence in the $N$-fold product of this basis.
Thus if $(e_n)$ is the original basis there exists a
permutation $\pi$ so that $(e_{\pi(n)})$ is equivalent to
the canonical basis of $\ell_F[s'_{1+[n-1/N]}].$
The argument of Proposition 5.4 again establishes that
$(e_{\pi(n)})$ is equivalent to $(e_n)$.
But now the new basis is equivalent to $(e_{Nn})$ which is
also equivalent to $(e_n).$\qed\enddemo

We will specialize to consider functions of the form
$F(t)\sim t^p|\log t|^{-a}$ where $p\ge 1$ and $a>0.$
More precisely let $g(\tau)=\min (1,\tau^{-1})$ and let
$F^{p,a}$ be the Orlicz function corresponding to
$\phi=p+ag$ i.e.  $ F^{p,a}(t)=t^{p+a}$ for $e^{-1}\le t\le
1$ and $F^{p,a}(t)=e^{-pa}t^p|\log t|^{-a }$ for $0<t\le
e^{-1}.$
These functions are convex and m-convex.

Now suppose $s_n\downarrow 0$.
For each $n\in\bold N$ let $N_n$ be the greatest index such
that $s_k\ge \exp(-2^n)$, and let $N_0=0.$
Let $E_n=\{N_{n-1}+1\le k\le N_n\}$ and $V_n=[e_k\colon k\in
E_n].$

\proclaim{Proposition 6.6}Suppose $1\le p<\infty$ and $a>0$
are fixed.
Let $F=F^{p,a}.$
Then if $0<s_n\le 1$ and $s_n\downarrow 0$, we have \newline
(1) $\ell_F[s_n]=\ell_p(V_n).$ \newline (2) There is a
constant $C$ depending only on $p,a,$ so that if $x\in V_n,$
then
$$
C^{-1}\|x\|_{\ell_{G_n}} \le \|x\|_{\ell_F[s_n]}\le
C\|x\|_{\ell_{G_n}}.
$$
\newline (3) $\ell_F[s_n]=\ell_p$ as a sequence space if and
only if there is a constant $K$ so that $N_n\le
\exp(K2^n).$\newline \endproclaim

\demo{Proof}(1) Notice that if $k\in E_n,$ then
$\sup_{\tau\ge |\log s_k|}\phi(\tau)\le g(2^{n-1})\le
p+a2^{-(n-1)}.$
Now by Lemma 6.1, if $x\in c_{00}$ we have
$$
\| (\|E_nx\|_{\ell_F[s_n]})\|_{\ell_{p+a2^{n-1}}}\le
\|x\|_{\ell_F[s_n]}\le \|(\|E_nx\|_{\ell_F[s_n]})\|_{\ell_p}
$$
and by the Simons criterion \cite{25} we obtain (1).

(2) If $k\in E_n,$ we have $F_{s_k}(t) \le
F_{\exp(-2^n)}(t)$ for $0\le t\le 1.$
Conversely
$$
F_{\exp(-2^n)}(t)\le \exp(\int_{\log
1/s_k}^{2^n}g(\tau)d\tau)F_{s_k}(t) \le eF_{s_k}(t).
$$

(3) If $x\in V_n$ then by Lemma 6.1 we have
$$
\|x\|_{\ell_{p+a2^{-(n-1)}}} \le \|x\|_{\ell_F[s_n]}\le
\|x\|_{\ell_p}
$$
so that if $N_n-N_{n-1}\le N_n\le \exp(K2^n)$ each
$(e_k)_{k\in E_n}$ is uniformly equivalent to the usual
basis of $\ell_p^{N_n-N_{n-1}}.$ Conversely note that
$$
\|\sum_{k=1}^{N_n}e_k\|_{\ell_F[s_m]} \le
N_n^{1/(p+a2^{-n})}
$$
so that the condition is also necessary.  \qed\enddemo

We now give a general criterion for checking permutative
equivalence of two bases in these special Orlicz modular
spaces.

\proclaim{Lemma 6.7}Suppose $1\le p<\infty$ and $a>0$ are
fixed and let $F=F^{p,a}.$
Suppose $0<s_n,s'_n\le 1$ and $s_n,s'_n\downarrow 0.$
Suppose the the canonical bases of $\ell_F[s_n]$ and
$\ell_F[s'_n]$ are permutatively equivalent.
Then there is a constant $K$ so that for every $n\ge 1$ and
$k\ge 2$ we have
$$
|\log s'_{n+k}|^{-1}\le |\log s_{n}|^{-1} + K(\log k)^{-1}
$$
and
$$
|\log s_{n+k}|^{-1} \le |\log s'_n|^{-1}+K(\log k)^{-1}.
$$
\endproclaim

\demo{Proof} Let us define
$$
\align D(n,k)&=\inf_{|A|=n+k}\sup_{B\subset A\atop |B|=k}
\|e_B\|_{\ell_F[s_m]}\\
D'(n,k)&=\inf_{|A|=n+k}\sup_{B\subset A\atop |B|=k}
\|e_B\|_{\ell_F[s'_m]}.  \endalign
$$
Then there is a constant $C$ so that for all $n,k$ we have
$C^{-1}D(n,k)\le D'(n,k)\le CD(n,k).$
Notice however that
$$
D(n,k) =\|\sum_{j=n+1}^{n+k}e_j\|_{\ell_F[s_n]}
$$
and hence if $\tau_j=\log 1/s_j,$
$$
k^{1/(p+ag(\tau_{n+k}))}\le D(n,k) \le k^{1/(p+ag(\tau_n))}
$$
and similarly for $D'(n,k).$
It follows that
$$
\frac{\log k}{p+ag(\tau_{n+k})} \le \log C + \frac{\log
k}{p+ag (\tau'_n)}
$$
and this combined with a similar inequality with roles
reversed gives the result.  \qed\enddemo

\proclaim{Proposition 6.8} Suppose $0<s_n\le 1$ and
$s_n\downarrow 0.$
The canonical basis of $\ell_F[s_n]$ is equivalent to its
square if and only if there exists $l\ge 1$ so that
$N_{n+l}+\exp(2^{n+l})\ge 2N_n$ for all $n\ge 1.$
\endproclaim

\demo{Proof} From the preceding lemma, we obtain that if the
canonical basis is equivalent to its square then,
$$
\log|\log s_{2n}| \le \log|\log s_n| + K(\log n)^{-1}
$$
for some constant $K.$
Now suppose $N_{n+l}+\exp(2^{n+l})\le 2N_n.$
Then $\log N_n \ge 2^{n+l}-\log 2\ge 2^{n+l-1}$ and hence
$$
|\log s_{N_n}|^{-1} \le |\log s_{2N_n}|^{-1}+2K2^{-n-l}\le
(1+2K)2^{-n-l}.
$$
Thus $ 2^{-n-1}\le (1+2K)2^{-n-l}$ so that $l\le
\log_2(2+4K).$
This implies the given criterion.

For the converse, notice that since the standard $\ell_p-$
basis is equivalent to some subsequence of the given basis,
the canonical basis is equivalent to the canonical basis of
a space $\ell_F[s'_n]$ where $N'_n=N_n+[\exp 2^n].$
It is then clear that for some fixed $l$ we have
$N'_{n+l}\le 2N'_n.$
This in turn implies that $|\log s'_{2n}|\le K|\log s'_{n}|$
for some constant $K$.
But then $F_{s'_{2n}}(t) \le K^aF_{s'_n}(t)$ for $0\le t\le
1$ whence the result.  \qed\enddemo

\proclaim{Theorem 6.9}Suppose $a>0$ and let $F(t)\sim t|\log
t|^{-a}$ for $t$ near zero.
Let $Z$ be a complemented subspace of $\ell_F$ with an
unconditional basis $(u_n).$
Suppose $Z$ is not isomorphic to $\ell_F.$
Then:\newline (1) If $(u_n)$ is equivalent to its square
then $Z$ has a unique unconditional basis.\newline
(2) If every complemented subspace of $Z$ with an
unconditional basis also has a unique unconditional basis
then $Z$ is isomorphic to $\ell_1.$ \endproclaim

\demo{Remark}By combining Propositions 6.6 and 6.8, it is
clear that we can find $(s_n)$ with $s_n\downarrow 0$, so
that the canonical basis is equivalent to its square, but
$\ell_F[s_n]$ is not isomorphic to $\ell_1.$
Thus Theorem 6.9 answers Problem 11.2 of \cite{3}
negatively.\enddemo

\demo{Proof} (1) has already been proved above; it is a
special case of Theorem 6.5.
For (2) we consider $F=F^{1,a}$ and a sequence $0<s_n\le 1$
with $s_n\downarrow 1.$ Let $N_n,E_n$ and $V_n$ be defined
as before and let $M_n=N_n-N_{n-1}.$
We may suppose, without loss of generality that
$s_k=\exp(-2^n)$ when $N_{n-1}\le k\le N_n,$ by applying
Proposition 6.6 (2).
Assume that every subsequence of the canonical basis of
$\ell_F[s_n]$ spans a space with a unique unconditional
basis.

Let $P_n=[\sqrt M_n].$
We use a result of \cite{15} that since the given basis of
each $V_n$ is symmetric there is an unconditional basis
$(u_k)_{k\in E_n}$ of each $V_n$ uniformly equivalent to the
direct sum of $M_n-P_n$ members of the given basis and $P_n$
constant coefficient vectors of length $P_n.$

Now if $\Cal N$ is any infinite subset of $\bold N$ we can
consider the basis $(u_k)_{k\in E_n,n\in\Cal N}$ of the
subspace $[e_k]_{k\in E_n,n\in \Cal N}.$
This is equivalent to the canonical basis of the space
$\ell_F[(s'_k)_{k\in E_n,n\in\Cal N}]$ where $s'_k=s_k$ for
$N_{n-1}+1\le k\le N_n-P_n$ and $s'_k=\rho_ns_k$ for
$N_n-P_n+1\le k\le N_n$ where $F(\rho_ns_k)=P_n^{-1}F(s_k).$
Clearly $\rho_n\le P_n^{-(1+a)^{-1}}.$

Now assume that $\ell_F[s_n]$ is not isomorphic to $\ell_1.$
Then $2^{-n}\log M_n$ is unbounded.
We then choose $\Cal N=\{n_1,n_2,\ldots\}$ inductively so
that $\rho_{n_j}\exp(2^{-n_j}) \ge \exp(2^{-n_{j+1}})$ for
$j=1,2,\ldots$ and that $2^{-n_j}\log M_j$ is unbounded.
Then the sequence $(s'_k)_{k\in E_n,n\in\Cal N}$ is already
in decreasing order and the corresponding basis is
equivalent to that for $(s_k)_{k\in E_n,n\in\Cal N}.$
Lemma 6.7 can now be used again to show that for some
constant $K$ we have that for $n\in\Cal N,$
$$
2^{-n}\le (|\log \rho_{n}|+2^{n})^{-1} + K(\log P_n)^{-1}).
$$
Since $|\log \rho_n|\ge (1+a)^{-1}\log P_n$ this implies
that
$$
\log P_n\le K'2^n
$$
for some $K'.$
Thus $\log M_n\le 3K'2^n$ for $n\in\Cal N$ and we have a
contradiction.\qed\enddemo

\demo{Remark}This theorem can be proved for wider range of
Orlicz functions.
Specifically a proof along the same lines can be given if
$\phi$ decreases monotonically, $\phi'$ is eventually
monotone increasing, $\phi(\tau)-1=O((\log \tau)^{-1})$ and
$\tau(\phi(\tau)-1)$ is eventually increasing.\enddemo

Finally let us notice it is also possible to give a
super-reflexive version.

\proclaim{Theorem 6.10}Let $F(t)\sim t^2|\log t|^{-1}$ for
$t$ near zero.
Let $Z$ be a complemented subspace of $\ell_F$ with an
unconditional basis $(u_n).$
Suppose $Z$ is not isomorphic to $\ell_F.$
Then:\newline (1) If $(u_n)$ is equivalent to its square
then $Z$ has a unique unconditional basis.\newline
(2) If every complemented subspace of $Z$ with an
unconditional basis also has a unique unconditional basis
then $Z$ is isomorphic to $\ell_2.$ \endproclaim

\demo{Proof}The proof of (2) is identical to the proof given
above.  For (1), we need a result analogous to Theorem 6.5.
An inspection of the proof reveals that it is only necessary
to show that every complemented unconditional basic sequence
is equivalent to a sequence of constant coefficient blocks.
It suffices to prove the same result in $\ell_F^*=\ell_G$
where $G(t)\sim t^2|\log t|.$
But every unconditional basic sequence in $\ell_G$ is
equivalent to sequence of constant coefficient blocks
\cite{6}.\qed\enddemo

\demo{Remark}In fact in the dual space $\ell_G$ the results
hold for any subspace with an unconditional basis (even if
uncomplemented).

It may also be shown that the theorem is valid for $F(t)\sim
t^2|\log t|^{-a}$ where $a\ge 1.$

This requires a complex interpolation technique which we
will not expound here.\enddemo

\enddocument
\bye


\vskip10pt \Refs

\ref\no 1\by S. Bellenot \paper Tsirelson superspaces and
$\ell_p$ \jour J. Functional Analysis \vol 69\yr 1986 \pages
207-228\endref

\ref\no 2\by B. Bollobas \book Combinatorics \publ Cambridge
University Press \yr 1986\endref


\ref\no 3 \by J. Bourgain, P.G.  Casazza, J. Lindenstrauss
and L. Tzafriri \book Banach spaces with a unique
unconditional basis, up to a permutation \bookinfo Memoirs
Amer.  Math.  Soc.  No. 322 \yr 1985\endref

\ref\no 4\by J. Bourgain, N.J.  Kalton and L. Tzafriri
\paper Geometry of finite-dimensional subspaces and
quotients of $L_p$ \jour Springer Lecture Notes 1376 \yr
1989 \pages 138-175\endref

\ref\no 5\by P.G.  Casazza, W.B.  Johnson and L.
Tzafriri\paper On Tsirelson's space \jour Israel J. Math.
\vol 47 \yr 1984 \pages 81-98\endref

\ref\no 6\by P.G.  Casazza and N.J.  Kalton \paper
Unconditional bases and unconditional finite-dimensional
decompositions in Banach spaces \jour Israel J. Math
\paperinfo to appear\endref

\ref\no 7 \by P.G.  Casazza, N.J.  Kalton and L. Tzafriri
\paper Uniqueness of unconditional and symmetric structures
in finite-dimensional spaces \jour Illinois J. Math.  \vol
34 \yr 1990 \pages 793-836\endref

\ref\no 8 \by P.G.  Casazza and T.J.  Schura \book
Tsirelson's space \bookinfo Springer Lecture Notes 1363
\publ Springer \yr 1989\endref


\ref\no 9 \by L.E.  Dor\paper On projections in $L_1$ \jour
Ann.  Math.  \vol 102 \yr 1975 \pages 463-474\endref

\ref\no 10\by I.S.  Edelstein and P. Wojtaszczyk \paper On
projections and unconditional bases in direct sums of Banach
spaces \jour Studia Math.  \vol 56 \yr 1976 \pages
263-276\endref


\ref\no 11\by W.T.  Gowers \paper A finite-dimensional
normed space with non-equivalent symmetric bases \jour
Israel J. Math.  \vol 87 \yr 1994 \pages 143-151\endref

\ref\no 12\by W.T.  Gowers \paper A solution to Banach's
hyperplane problem\jour Bull.  London Math.  Soc.  \vol 26
\yr 1994 \pages 523-530 \endref

\ref\no 13\by W.T.  Gowers and B. Maurey \paper Banach
spaces with small spaces of operators \paperinfo
preprint\endref


\ref\no 14\by F.L.  Hernandez and N.J.  Kalton \paper
Subspaces of rearrangement-invariant spaces\paperinfo to
appear\endref

\ref\no 15 \by W.B.  Johnson, B. Maurey, G. Schechtman and
L. Tzafriri \book Symmetric structures in Banach spaces
\bookinfo Mem. Amer.  Math.  Soc.  No. 217 \yr 1979\endref







\ref\no 16\by N.J.  Kalton\book Lattice structures on Banach
spaces \bookinfo Mem. Amer.  Math.  Soc.  No. 493 \yr
1993\endref

\ref\no 17 \by G. K\"othe and O. Toeplitz \paper Lineare
Raume mit unendlich vielen Koordinaten und Ringen
unendlicher Matrizen \jour J. Reine Angew.  Math.  \vol 171
\yr 1934 \pages 193-226 \endref

\ref\no 18 \by J. Lindenstrauss and A. Pe\l czynski \paper
Absolutely summing operators in $\Cal L_p$-spaces and their
applications \jour Studia Math.  \vol 29 \yr 1968 \pages
315-349 \endref

\ref\no 19 \by J. Lindenstrauss and L. Tzafriri \book
Classical Banach spaces I, Sequence spaces \publ
Spring\-er
Verlag\publaddr Berlin, Heidelberg, New York \yr
1977\endref


\ref\no 20 \by J. Lindenstrauss and L. Tzafriri \book
Classical Banach spaces II, Function spaces \publ Spring\-er
Verlag\publaddr Berlin, Heidelberg, New York \yr 1979\endref



\ref\no 21 \by J. Lindenstrauss and M. Zippin \paper Banach
spaces with a unique unconditional basis \jour J. Functional
Analysis \vol 3 \yr 1969 \pages 115-125\endref

\ref\no 22 \by B.S.  Mityagin \paper Equivalence of bases in
Hilbert scales (in Russian) \jour Studia Math.  \vol 37 \yr
1970 \pages 111-137\endref

\ref\no 23\by G. Pisier \book The volume of convex bodies
and geometry of Banach spaces \bookinfo Cambridge Tracts
94\publ Cambridge University Press\yr 1989\endref


\ref\no 24\by C. Sch\"utt \paper On the uniqueness of
symmetric bases in finite dimensional Banach spaces \jour
Israel J. Math.  \vol 40 \yr 1981 \pages 97-117 \endref

\ref\no 25 \by S. Simons \paper The sequence spaces
$\ell(p_{\nu})$ and $m(p_{\nu})$ \jour Proc.  London Math.
Soc.  \vol (3) 15 \yr 1965 \pages 422-436\endref


\ref\no 26\by P. Wojtaszczyk \paper Uniqueness of
unconditional bases in quasi-Banach spaces with applications
to Hardy spaces, II \jour Israel J. Math \paperinfo to
appear \endref

\ref\no 27 \by M. Wojtowicz \paper On Cantor-Bernstein type
theorems in Riesz spaces \jour Indag.  Math.  \vol 91\yr
1988 \pages 93-100\endref


\endRefs
\bye